
\def\LaTeX{%
  \let\Begin\begin
  \let\End\end
  \let\salta\relax
  \let\finqui\relax
  \let\futuro\relax}

\def\UK{\def\our{our}\let\sz s}
\def\USA{\def\our{or}\let\sz z}



\LaTeX

\USA


\salta

\documentclass[twoside,12pt]{article}
\setlength{\textheight}{24cm}
\setlength{\textwidth}{16cm}
\setlength{\oddsidemargin}{2mm}
\setlength{\evensidemargin}{2mm}
\setlength{\topmargin}{-15mm}
\parskip2mm


\usepackage[usenames,dvipsnames]{color}
\usepackage{amsmath}
\usepackage{amsthm}
\usepackage{amssymb,bbm}
\usepackage[mathcal]{euscript}

\usepackage{hyperref}
\usepackage{enumitem}

%
%

%
 
\definecolor{rosso}{rgb}{0.85,0,0}
\definecolor{electricviolet}{rgb}{0.56, 0.0, 1.0}
\definecolor{pinegreen}{rgb}{0.0, 0.47, 0.44}

\def\juerg #1{{#1}}
\def\an #1{{\color{red}#1}}
\def\an #1{{#1}}
\def\pier #1{{\color{red}#1}} 
\def\pier #1{{#1}} 
\def\pcol #1{{\color{red}#1}}
\def\pcol #1{{#1}}
\def\pico #1{{\color{blue} #1}}
\def\pico #1{{#1}}




\bibliographystyle{plain}


%
\newtheorem{theorem}{Theorem}[section]
\newtheorem{remark}[theorem]{Remark}
\newtheorem{corollary}[theorem]{Corollary}

\finqui
\let\non\nonumber




\def\step #1 \par{\medskip\noindent{\bf #1.}\quad}
\def\jstep #1: \par {\vspace{2mm}\noindent\underline{\sc #1 :}\par\nobreak\vspace{1mm}\noindent}

\def\Lip{Lip\-schitz}
\def\Holder{H\"older}
\def\Frechet{Fr\'echet}
\def\lhs{left-hand side}
\def\rhs{right-hand side}


\def\multibold #1{\def\arg{#1}%
  \ifx\arg\pto \let\next\relax
  \else
  \def\next{\expandafter
    \def\csname #1#1#1\endcsname{{\bf #1}}%
    \multibold}%
  \fi \next}

\def\pto{.}

\def\multical #1{\def\arg{#1}%
  \ifx\arg\pto \let\next\relax
  \else
  \def\next{\expandafter
    \def\csname cal#1\endcsname{{\cal #1}}%
    \multical}%
  \fi \next}


\def\multimathop #1 {\def\arg{#1}%
  \ifx\arg\pto \let\next\relax
  \else
  \def\next{\expandafter
    \def\csname #1\endcsname{\mathop{\rm #1}\nolimits}%
    \multimathop}%
  \fi \next}

\multibold
qwertyuiopasdfghjklzxcvbnmQWERTYUIOPASDFGHJKLZXCVBNM.

\multical
QWERTYUIOPASDFGHJKLZXCVBNM.

\multimathop
diag dist div dom mean meas sign supp .


\def\Accorpa #1#2 #3 {\gdef #1{\eqref{#2}--\eqref{#3}}%
  \wlog{}\wlog{\string #1 -> #2 - #3}\wlog{}}


\def\<#1>{\mathopen\langle #1\mathclose\rangle}
\def\norma #1{\mathopen \| #1\mathclose \|}

\def\I2 #1{\int_{Q_t}|{#1}|^2}
\def\IT2 #1{\int_{Q_t^T}|{#1}|^2}
\def\IO2 #1{\norma{{#1(t)}}^2}
\def\IOT2 #1{\norma{{#1(T)}}^2}
\def\ov #1{{\overline{#1}}}

\def\iot {\int_0^t}
\def\ioT {\int_0^T}
\def\intQt{\int_{Q_t}}
\def\intQ{\int_Q}
\def\iO{\int_\Omega}

\def\Qtt{\int_{Q_t^T}}

\def\dt{\partial_t}
\def\ddt{\partial_{tt}}
\def\dn{\partial_{\bf n}}
\def\checkmmode #1{\relax\ifmmode\hbox{#1}\else{#1}\fi}

\def\bh{{\bf h}}

\def\erre{{\mathbb{R}}}

\def\II{{\mathbb{I}}}
\def\Vp{{V^*}}

\def\J{{\cal J}}

\def\Jred{{\J}_{\rm red}}

\def\S{{\cal S}}
\def\X{{\cal X}}
\def\Y{{\cal Y}}
\def\UU{{\cal U}}

\def\CP{{${\bf CP}$}}

\def\aeO{\checkmmode{a.e.\ in~$\Omega$}}
\def\aeQ{\checkmmode{a.e.\ in~$Q$}}




\def\genspazio #1#2#3#4#5{#1^{#2}(#5,#4;#3)}
\def\spazio #1#2#3{\genspazio {#1}{#2}{#3}T0}

\def\L {\spazio L}
\def\H {\spazio H}
\def\W {\spazio W}

\def\C #1#2{C^{#1}([0,T];#2)}



\def\Lx #1{L^{#1}(\Omega)}
\def\Hx #1{H^{#1}(\Omega)}

\def\Ldue{\Lx 2}

\def\Huno{\Hx 1}
\def\Hdue{\Hx 2}



\let\eps\varepsilon
\let\vp\varphi

\def\a{\alpha}	
\def\b{\beta}	
\def\d{\delta} 
\def\et{\eta}  
\def\th{\theta}

\def\ph{\varphi}

\def\ps{\psi}     

\def\cd{C_{\d}}

\def\thc{\theta_c}

\let\TeXchi\chi                         
\newbox\chibox
\setbox0 \hbox{\mathsurround0pt $\TeXchi$}
\setbox\chibox \hbox{\raise\dp0 \box 0 }
\def\chi{\copy\chibox}



\let\hat\widehat

\def\Uad{{\cal U}_{\rm ad}}
\def\UR{{\cal U}_{R}}

\def\bvp{\overline\varphi}

\def\bph{{\ov \ph}}

\def\ds{{\rm d}s}
\def\bu{{\ov u}}   
\def\bw{{\ov w}}  
\def\bdtw{{\dt \ov w}}    
\def\bv{{\ov v_0}}   
\def\opt{{(\bu,\bv)}}
\def\bwh{{\bw^{\bh}}}
\def\bvph{{\bvp^{\bh}}}

\usepackage{amsmath}
\DeclareFontFamily{U}{mathc}{}
\DeclareFontShape{U}{mathc}{m}{it}%
{<->s*[1.03] mathc10}{}

\DeclareMathAlphabet{\mathscr}{U}{mathc}{m}{it}
\Begin{document}


%
\title{Analysis and optimal control theory for a \\ \pcol{phase field model of Caginalp type \\ with thermal memory}}
\author{}
\date{}
\maketitle
\begin{center}	
\vskip-1.5cm
{\large\sc Pierluigi Colli$^{(1)}$}\\
{\normalsize e-mail: {\tt \href{mailto:pierluigi.colli@unipv.it}{pierluigi.colli@unipv.it}}}\\[0.25cm]
{\large\sc Andrea Signori$^{(1)}$}\\
{\normalsize e-mail: {\tt \href{mailto:andrea.signori01@unipv.it}{andrea.signori01@unipv.it}}}\\[0.25cm]

{\large\sc J\"urgen Sprekels$^{(2)}$}\\
{\normalsize e-mail: {\tt \href{mailto:juergen.sprekels@wias-berlin.de}{juergen.sprekels@wias-berlin.de}}}\\[.5cm]
$^{(1)}$
{\small Dipartimento di Matematica ``F. Casorati''}\\
{\small Universit\`a di Pavia}\\
{\small via Ferrata 5, I-27100 Pavia, Italy}\\[.3cm] 
$^{(2)}$
{\small Department of Mathematics}\\
{\small Humboldt-Universit\"at zu Berlin}\\
{\small Unter den Linden 6, D-10099 Berlin, Germany}\\[2mm]
{\small and}\\[2mm]
{\small Weierstrass Institute for Applied Analysis and Stochastics}\\
{\small Mohrenstrasse 39, D-10117 Berlin, Germany}\\[10mm]
\end{center}

%

\Begin{abstract}
\noindent A nonlinear extension of the Caginalp phase field system is considered that takes thermal memory into account.
The resulting model, which is a first-order approximation of a thermodynamically consistent system, is 
inspired by the theories developed by Green and Naghdi.  Two equations, resulting from phase dynamics and the universal  balance law for internal energy, are written in terms of the phase variable 
\juerg{(representing a non-conserved order parameter)} and the so-called thermal displacement, i.e., a primitive with respect to time of temperature. Existence and continuous dependence results are shown for weak and strong solutions to the corresponding initial-boundary value problem. \an{Then, an} optimal control problem is investigated for a suitable cost functional, in which two data act as controls, namely, the distributed heat source and the initial temperature. 
\Frechet\ differentiability between suitable Banach spaces is shown for the control-to-state operator, and meaningful 
first-order necessary optimality conditions are derived in terms of variational inequalities involving the adjoint variables.  \an{Eventually, c}haracterizations of the optimal controls are given. 

\End{abstract}
\salta
\pagestyle{myheadings}
\newcommand\testopari{\sc Colli--Signori--Sprekels}
\newcommand\testodispari{\sc Analysis and optimal control of a Caginalp type system
}
\markboth{\testopari}{\testodispari}
\finqui
%
\section{Introduction}
This paper is concerned with a phase field model for a non-isothermal phase transition
\juerg{with non-conserved order parameter} 
describing the evolution \juerg{in a container} in terms of two physical variables. Well-posedness issues for weak and strong solutions and optimal control problems are investigated in detail. At first, we 	
introduce the system of partial differential equations and related conditions.

\subsection{The initial and boundary value problem}
\juerg{We assume that the phase transformation takes place in a fixed container}
$\Omega\subset\erre^d$, $d \in \{2,3\}$, \juerg{which is} an open \pier{and bounded domain}
with smooth boundary $\Gamma:= \partial \Omega$. For a positive fixed final time horizon $T$, we set, 
\an{
\begin{align*}
	Q_t:=\Omega\times (0,t),\   \text{$\, 0< t \leq T$, } \quad
	Q:= Q_T,\quad \Sigma:=\Gamma\times (0,T).
\end{align*}}
Then the \juerg{model under study} reads as
\begin{alignat}{2}
\label{sys1}
&\dt\ph-\Delta\ph + \gamma(\ph) + \tfrac 2 {\thc} \pi(\ph) - \tfrac 1{\thc^2}{\dt w}\,\pi(\ph) \ni 0 \qquad&&\mbox{in }\,Q,\\
\label{sys2}
&\ddt w -\alpha \Delta (\dt w)-\beta \Delta w + \pi (\ph) \dt \ph =\an{u} \qquad&&\mbox{in }\,Q,\\
\label{sys3}
&\dn \ph=\dn (\alpha \dt w + \beta w)=0 \qquad&&\mbox{on }\,\Sigma,\\
\label{sys4}
&\ph(0)=\ph_0,\quad w(0)=w_0,\quad \dt w(0)=\an{v_0}\,\qquad &&\mbox{in }\,\Omega.
\end{alignat}
\Accorpa\Sys {sys1} {sys4}
The primary variables of the system are $\ph$, the  order parameter
\juerg{of the phase transition}, and $w$, the so-called {\em thermal displacement} or 
{\em freezing index}. 
The latter is directly connected to the \pier{absolute} temperature $\th$ of the system through the relation
\begin{align}
	\label{thermal_disp}
	w (\cdot , t)  = w_0 + \iot \th(\cdot, s) \,\ds, \quad t \in[0,T].
\end{align}
Moreover, $\a$ and $\b$ stand for prescribed \juerg{positive coefficients that are} related to the heat flux, $\thc$ for a  (positive) critical temperature, and $\an{u}$ for a distributed heat source.
Besides, the nonlinearities $\gamma: \erre \to 2^\erre$ and $\pi: \erre \to \erre$ indicate, in this order, a maximal monotone graph and a \Lip\ continuous function.
\juerg{Finally}, the symbol $\dn$ represents the outward normal derivative on $\Gamma$, whereas $\ph_0,w_0,$ and $\an{v_0}$ stand for some prescribed initial values.

Notice that the inclusion in \eqref{sys1} 
\juerg{is of Allen--Cahn type and is suited for the case of non-conserved order parameters (while the
case of a conserved order parameter would require a Cahn--Hilliard structure). The inclusion
originates from} the possibly multivalued nature of the graph~$\gamma$.
Typically, the maximal monotone graph $\gamma$ is obtained as the subdifferential of a 
convex and lower semicontinuous function $\hat \gamma : \erre \to [0,+\infty]$, and 
well-known examples are given by the regular, logarithmic, and double obstacle 
potentials, \an{defined, in the order,~by}
\begin{align}
  \label{Freg}
  \hat\gamma_{reg}(r)&=\frac {r^4}4\,, \quad r\in\erre\,,\\[2mm]
  \label{Flog}
  \hat\gamma_{log}(r)&=
 \begin{cases}
		\frac{\kappa}{2}[(1+r)\ln(1+r)+(1-r)\ln(1-r)],  &\hbox{if $r \in (-1,1)$},
		\\[1mm]
		\juerg{\kappa \ln(2)} ,  &\hbox{if $r \in \{-1,1\}$},
		\\[.5mm]
		+\infty,  &\hbox{otherwise},
		\end{cases}	
  \\[2mm]
  \label{F2ob}
  \hat\gamma_{dob}(r)&=I_{[-1,1]}(r),
\end{align}
with a positive constant $\kappa$, where, for every subset $A \subset \erre$, $I_{A}(\cdot)$ stands for the indicator function of \pico{$A$ and is} specified by 
\begin{align*}
	I_{A}(r):= 
	\begin{cases}
	0 \quad & \text{if $r \in A$},
	\\
	+ \infty & \text{otherwise}.
	\end{cases}
\end{align*}
Let us point out that the \an{inclusion} \eqref{sys1} simply reduces to an equality 
in the case of \eqref{Freg} and of \eqref{Flog} for $-1< \ph < 1$, since the regularity 
of $\hat \gamma$ ensures $\gamma$ to be single valued.

Next, we present a \an{possible} physical derivation of the system in \eqref{sys1}--\eqref{sys4},
trying to \juerg{meet the requirement of thermodynamic} consistency as much as possible. 
On the other hand, 
different approaches may be appealed and, in particular, we quote 
\cite{BR17, BS, Cag, CM, KS04, MQ10, PG09} as related references.


\subsection{Thermodynamic derivation and modeling considerations}

We start from the local \juerg{specific} Helmholtz free energy, acting on the \juerg{absolute temperature} $\th >0$ and the \juerg{dimensionless} order parameter $\ph$.  \juerg{With physical constants 
$\beta_1,\beta_2,\beta_3$, the specific local free energy $F$ is assumed in the form}
\begin{equation} \label{free}
  F(\th,\ph)=
 {c_V}\th (1-\ln (\th/\th_1)) + \juerg{\beta_1}\hat \pi (\ph) +\juerg{\beta_2}\th \, \hat \gamma (\ph) 
 +\juerg{\frac{\beta_3}2}\,\th\,|\nabla\varphi|^2, 
\end{equation}
\juerg{
where $c_V >0 $ denotes the specific heat (assumed constant), 
$\th_1>0 $ is some fixed reference temperature, $\hat \gamma (\ph)$ has been introduced
above,} and the real-valued function $\hat \pi$ stands for a primitive of $\pi$.
The last summand in \eqref{free} is a contribution that accounts for nearest-neighbor interactions. 

\juerg{By virtue of the general relations between the thermodynamic potentials, the expressions for local 
 specific entropy~$S$ and local specific internal energy~$E$ are then given by}
\begin{align}
 S(\th,\ph) &= - \partial_\th  F(\th,\ph) = {c_V}\ln (\th/\juerg{\th_1})  
 -\juerg{\beta_2}\,\hat \gamma (\ph)  - \frac{\juerg{\beta_3}}{2}|\nabla\varphi|^2, 
\label{umbazero}\\
 E(\th,\ph) &= F(\th,\ph) + \th  S(\th,\ph) = c_V \th + \juerg{\beta_1}\hat \pi (\ph).
\label{umbauno}
\end{align} 

Now, we come to the evolution laws. \juerg{As always, the universal balance law of internal
energy must be obeyed. Under the assumption that velocity effects may be discarded, it has the general form}
\begin{align}
  &\juerg{\rho}\,\partial_t E(\th,\ph) +\div \mathbf q = \juerg{\rho}\hskip1pt \an{u},
  \label{umbadue}
\end{align}
\juerg{where \,$\mathbf q$\,  denotes the heat flux,  $\rho$\, is the mass density and \,$\rho \hskip1pt\an{u}$\, stands 
for the possible presence of distributed heat sources/sinks. Here, we consider the case when $\,\rho\,$ varies only little during the phase transition and can be assumed constant. }

Usually the Fourier law is assumed for $\bf q$, i.e., 
\begin{equation}
  \label{umbatre}
\mathbf q=-\kappa_V\nabla\th,
\end{equation}
\juerg{where $\kappa_V$ is the (positive) heat conductivity} coefficient, 
together with the no-flux condition $\mathbf q\cdot \mathbf n=0$ on the boundary.

In the present paper, we adopt a different approach for $\bf q$,
the Fourier law~\eqref{umbatre} being generalized in the light of 
the works by Green and Naghdi~\cite{GN91,GN92,GN93} and Podio-Guidugli~\cite{PG09}.
Indeed, these authors introduced a different approach for the study of 
heat conduction theory \juerg{that leads to the notion of {\em thermal displacement}}. 
We recall \eqref{thermal_disp} and note that there $w_0$ represents 
a given datum at the (initial) reference time. This datum accounts for a possible previous
thermal history of the phenomenon.
Making use of this new variable~$w$, 
Green and Naghdi proposed three theories for heat transmission labeled as
type I--III. Let us now employ the symbols $\alpha$ and $\beta$ for the coefficients
\juerg{which are assumed constant and positive}. 
Type I theory, after suitable linearization, brings us back to the standard Fourier~law
\begin{equation}
  \mathbf q=-\alpha\nabla (\dt w) \qquad\text{(type I),} 
  \label{eq:5}
  \end{equation}
while linearized versions of type II and III yield the following
heat-conduction laws:
\begin{alignat}{2}
  &\mathbf q=-\beta \nabla w \qquad &&\text{(type II),} \label{eq:6}
    \\&
  \mathbf q=-\alpha \nabla (\dt w )- \beta \nabla w \qquad &&\text{(type III).}
  \label{eq:7}
\end{alignat}
We point out that the thermal displacement $w$ is useful 
to describe type II and III laws, whereas the type I law can be
stated in terms of the temperature $\th=\partial_t w$ alone.

This paper is concerned with the general type III theory. In fact, 
in view of \eqref{umbadue} and \eqref{eq:7}, we infer \an{that}
\begin{equation}
\juerg{\rho}\,\bigl(c_V w_{tt} + \juerg{\beta_1}\pi (\ph) \dt \ph\bigr) -\alpha \Delta (\dt w)-\beta \Delta w =
\juerg{\rho}\hskip1pt\an{u}. 
\label{eq:8}
\end {equation}
\juerg{Observe that the no-flux condition $\mathbf q\cdot\mathbf n=0$ then gives rise to the second
boundary condition in \eqref{sys3}.}

\juerg{It remains to derive the equation governing the evolution of the order parameter. 
To this end, we introduce the total entropy functional, which at any
fixed time instant $t\in[0,T]$ is given by the expression}
$$ \juerg{{\mathscr S}[\th(t),\ph(t)]= \int_\Omega \rho\, S(\th(t),\ph(t))\,,} 
$$
\juerg{with the usual notation $\th(t)=\th(\cdot,t), \ph(t)=\ph(\cdot,t)$.}

\juerg{For the dynamics of the order parameter, we postulate that it
runs at each time instant $t\in(0,T]$ in a direction as to maximize total entropy subject
to the constraint that the balance law \eqref{umbadue} of internal energy be satisfied.
To this end, observe that integration of \eqref{umbadue} over $\Omega\times[0,t]$, using 
\eqref{umbauno} and the
no-flux boundary condition for $\mathbf q$, yields the identity}
$$ 
\juerg{0\,=\,\iO\rho\bigl(c_V \an{\theta(t)}-c_V\theta_0+\beta_1
\widehat \pi(\an{\ph(t)})-\beta_1\widehat\pi(\ph_0)-\an{R(t)}
\bigr)}, 
$$
\juerg{where \an{we again use the notation $R(t)=R(\cdot,t)$, and} $\,R(x,t)\an{:}=\int_0^t \an{u}(x,s)\,{\rm d}s$\an{, $x \in \Omega$}. 
We now consider the augmented entropy functional }
\begin{align*}
{\mathscr S}_\lambda[\th(t),\ph(t)] :=\, &{\mathscr S}[\th(t),\ph(t)] + 
\juerg{\rho\iO\lambda(\cdot,t)\,\bigl(c_V\theta(t)-c_V\theta_0+\beta_1
\widehat\pi(\ph(t))-\beta_1\widehat\pi(\ph_0)-\an{R(t)}
\bigr)}\\
=\,&\juerg{\rho\iO\bigl[c_V\ln(\theta(t)/\theta_1)-\beta_2\,\widehat\gamma(\ph(t))-\,\frac{\beta_3}2\,
|\nabla\ph(t)|^2 }\\
&\juerg{\hspace*{10mm} +\,\lambda(\cdot,t)\,\bigl(c_V\theta(t)-c_V\theta_0+\beta_1
\widehat\pi(\ph(t))-\beta_1\widehat\pi(\ph_0)-\an{R(t)} \bigr)\bigr], }
\end{align*}
where \pico{$\lambda(t) = \juerg{\lambda(x,t)} $, $x\in\Omega$,} plays the role of a \an{Lagrange} multiplier. The search for critical points leads to the Euler--Lagrange equations obtained by taking the variational derivatives of $ {\mathscr S}_\lambda$ with respect to $\ph $ and $\th$, namely, 
\begin{align}
&\delta_\ph {\mathscr S}_\lambda\juerg{[\th(t),\ph(t)] = \,\rho\bigl[ -\beta_2\gamma (\ph(t))+\beta_3\Delta\ph(t) 
 + \lambda(t)\,\beta_1 \pi (\ph(t))\bigr]} \ni 0 , \non \\
&\delta_\th {\mathscr S}_\lambda\juerg{[\th(t),\ph(t)] = \rho\bigl[{c_V}/{\th(t)} +\lambda(t)\, c_V\bigr]} =0 . \non
\end{align} 
Then, from the second relation we can identify $\lambda $ as $ - 1/\th$, while we 
postulate that the evolution of $\ph$ runs 
in the direction of $\delta_\ph {\mathscr S}_\lambda $ 
at a rate which is proportional to it. More precisely, 
we assume that the evolution of $\ph$ is governed by the equation\pico{%
$$a_V(\theta, \ph) \partial_t \varphi \,= \delta_\ph {\mathscr S}_\lambda(\theta, \ph),$$
that corresponds to   
\begin{align}
  a_V(\theta, \ph) \partial_t \varphi \, 
    = \,\juerg{\rho\bigl[- \pico{(\beta_1/\theta)} \pi (\ph) - \beta_2\,\gamma (\ph) + \beta_3\,\Delta \varphi \bigr],}
  \label{eq:1}
\end{align}
}%
\juerg{where $\,a_V\,$ is a positive coefficient (assumed constant).}

\juerg{At this point, we simplify the exposition by generally assuming in the following that the numerical values of all of the physical constants $c_V,\rho,\beta_1,\beta_2,\beta_3, a_V$ equal unity, while their physical dimensions
will be kept active so that they still match. This will have no bearing on the
subsequent mathematical analysis and should not lead to any confusion. However, 
in a practical application of the model with real 
physical data, this would have to be accounted for. Under these premises, 
the balance of internal energy \eqref{eq:8} takes the form \eqref{sys2}, 
and \eqref{eq:1} becomes}
\begin{align}
  \partial_t \varphi   - \Delta \varphi  + \juerg{\gamma (\ph)} + \frac1 \theta \pi (\ph) \pico{{}\ni{}} 0. \label{eq:2}
\end{align}
From \eqref{eq:2} we arrive at \eqref{sys1} with the help of \eqref{thermal_disp} and of the first-order 
approximation 
$$ \frac 1 \theta \approx \frac 1 \thc - \frac1 {\thc^2} (\th - \thc )$$
about the critical temperature $\thc$.

Initial conditions  for $\ph, \, w, \, \dt w$ are prescribed in \eqref{sys4} to complete the initial boundary value problem.


\subsection{Comments and results}

The full set \an{of equations \Sys\ turns} out to be a variation of the Caginalp phase field model~\cite{Cag}. Some 
mathematical discussion of a simpler problem for \eqref{sys1}--\eqref{sys2} has already been 
given in~\cite{MQ10}. The papers~\cite{CC12,CC13} dealt with well-posedness
issues and asymptotic analyses with respect to the positive coefficients $ \alpha, \, 
\beta$ as one of them \juerg{approaches zero}. Other concerned results for this class of systems 
\an{may be found} in \cite{CDM, CGM13}. \an{Finally, let us notice that sliding} mode control problems 
were investigated in \cite{CM}. 

The existence of a weak solution for \eqref{sys1}--\eqref{sys4} and its continuous dependence
with respect to data are \juerg{for the first time} examined in the present paper, under very general assumptions on the convex function $\hat \gamma$. Then, the regularity issue for obtaining  strong solutions of the system is analyzed and an improved continuous
dependence estimate is proved in a restricted framework for $\hat \gamma$ \juerg{that still allows for} the cases 
\eqref{Freg} and \eqref{Flog} of regular and logarithmic potentials. However,
the point of emphasis for this paper is the study of the optimal control problem,
whose precise formulations is given at the beginning of Section~\ref{SEC:OC} (cf.~\eqref{cost}--\eqref{Uad}). A tracking-type functional has to be minimized  
with respect to the variation of the distributed heat source $\an{u}$ in 
\eqref{sys2} and of the initial value $\an{v_0} $ 
for the temperature $\dt w$. Indeed, both these data are taken as controls, 
and the existence of optimal controls is investigated along with first-order necessary optimality conditions. More specifically, the linearized problem is introduced, and it is shown that the control-to-state mapping is Fr\'echet differentiable between suitable spaces. The optimal controls are \an{eventually} characterized in terms of variational inequalities for the associated adjoint variables.

About optimal control problems for phase field systems, in particular of Caginalp type, 
we can quote the pioneering work \cite{HJ92}; one may also see the specific sections in the monograph~\cite{Trol}. For other contributions, we mention the article\cite{LS}, dedicated to a thermodynamically consistent version of the phase field system described above, and the more recent papers\cite{CGMR1} and \cite{CGS26}, where the interested reader can find a list of related references.

\subsection{Preliminaries }
Let us set the notation we are going to employ throughout the paper.
Given a Banach space $X$, we denote by $\norma{\cdot}_X$ the corresponding norm, by $X^*$ its topological dual space, and by $\< \cdot , \cdot >_X$ the related duality pairing between $X^*$ and $X$. 
The standard Lebesgue and Sobolev spaces defined on $\Omega$, for every $1 \leq p \leq \infty$ and $k \geq 0$, are denoted by $L^p(\Omega)$ and $W^{k,p}(\Omega)$, and the associated norms by $\norma{\cdot}_{L^p(\Omega)}=\norma{\cdot}_{p}$ and $\norma{\cdot}_{W^{k,p}(\Omega)}$, respectively. 
For the special case $p = 2$, these become Hilbert spaces, and we denote by $\norma{\cdot}=\norma{\cdot}_2$ the norm of $\Lx2$ and employ the usual notation $H^k(\Omega):= W^{k,2}(\Omega)$. 

For convenience, we also \juerg{introduce the notation}
\begin{align}
  & H := \Ldue \,, \quad  
  V := \Huno\,,   \quad
  W := \{v\in\Hdue: \ \dn v=0 \,\mbox{ on $\,\Gamma$}\}.
  \label{def:HVW}
\end{align}
Besides, for Banach spaces $X$ and $Y$, we introduce the linear space
$X\cap Y$, which becomes a Banach space when equipped with its natural norm $\,\|v\|_{X \cap Y}:=
\|v\|_X + \|v\|_Y$, for $v\in X\cap Y$. 
To conclude, \pier{for normed spaces $\,X\,$ and $\,v\in L^1(0,T;X)$, we set}
\begin{align}
	(1 * v)(t):=\int_0^t v(s)\,\ds,\quad \hbox{$t \in[0,T]$},
	\label{intr1}
\end{align}
\pier{and also introduce the notation
\begin{align}
	(1 \circledast v)(t):=\int_t^T v(s)\,\ds\an{,\quad \hbox{$t \in[0,T]$}.}
	\label{intr2}	
\end{align}}%
Throughout the paper, we employ the following convention: the capital-case symbol $C$ is used to denote every constant
that depends only on the structural data of the problem such as
$T$, $\Omega$, $\a$, $\b$, $\thc$, the shape of the 
nonlinearities, and the norms of the involved functions. For this reason, its meaning may vary from line to line
\juerg{and even within formulas}.
Moreover, when a positive constant $\delta$ enters the computation, the related symbol $\cd$ denotes constants
that depend on $\d$ in addition.

\subsection{Plan of the paper}
The rest of the work is organized in the following way. Section \ref{SEC:WP} is devoted to 
the mathematical analysis of system \Sys. We prove the existence and uniqueness of a weak solution in a very general framework \juerg{that} includes singular and nonregular potentials like the double obstacle one. We then show that in the case of regular and logarithmic potentials, under natural assumptions for the initial data, the system admits a unique strong solution and that the phase variable enjoys the so-called separation property. This \an{latter} is of major importance for the mathematical analysis of phase field models involving singular potentials \an{as it guarantees that the singularity of the potential $\gamma$ is no longer an obstacle for the mathematical analysis. In fact, it ensures the phase field variable $\ph$ to range in some interval in which the potential is smooth.}
Next, in Section \ref{SEC:OC}, by the results \pier{shown} in Section \ref{SEC:WP}, we discuss a nontrivial application to optimal control, where we seek optimal controls in the form of a distributed heat source and an initial temperature. \an{The} existence of an optimal strategy as well as first-order necessary 
optimality conditions \an{are addressed}.

\section{Analysis of the system}
\label{SEC:WP}
\setcounter{equation}{0}

The following assumptions will be in order throughout this paper.
\begin{enumerate}[label={\bf A\arabic{*}}, ref={\bf A\arabic{*}}]
\item \label{ass:1:constants}
	$\a,\b,$ and $\thc$ are positive constants.
\item \label{ass:2:gamma}
	$\hat \gamma:\erre\to [0,+\infty]$ is convex and lower semicontinuous with
	$\hat \gamma(0)=0$, so that $\gamma := \partial\hat \gamma$ is a maximal monotone graph with $\gamma (0) \ni 0$.
	Moreover, we denote the effective domain of $\gamma$ by $\dom(\gamma).$
\item \label{ass:3:pi}
		\pier{$\pi: \erre \to \erre$ is a \Lip\ continuous function.
		Let $\hat \pi \in C^1(\erre)$ denote a primitive of $\pi$, i.e., $\pi (r) = \hat \pi ^{\, \prime} (r) $ for every $r\in \erre$.} 
\end{enumerate}

The first result concerns the existence of \pier{weak solutions}. 
\begin{theorem}
\label{THM:EX:WEAK}
Assume \an{that \ref{ass:1:constants}--\ref{ass:3:pi} hold}.
Moreover, let the initial data fulfill
\begin{align}
	\label{ass:exweak:initialdata}
	\ph_0 \in V, \quad \hat \gamma(\ph_0)\in\Lx1, \quad w_0 \in V, \quad \an{v_0} \in H,
\end{align}
and\pier{, for the heat source,} \an{suppose that}
\begin{align}
	\label{ass:regf:weak}
	\an{u} \in \L2 H. 
\end{align}
Then there exists a weak solution $(\ph, w, \xi)$ to \pier{the} system \Sys\ in the sense that 
\begin{align*}
	\ph &\in \H1 H \cap \L\infty V \cap \L2 W,
	\\ 
	\xi & \in	\L2 H,
	\quad \ph \in \dom(\gamma) \quad \text{and} \quad \xi \in \gamma(\ph) \,\, a.e. \,\, \text{in} \,\, Q,
	\\ 
	w & \in \H2 {\Vp}\cap\W{1,\infty} H \cap \H1 V,
\end{align*}
and \pier{that the variational equalities}
\begin{align}
	\label{wf:1}
	& \iO \dt\ph \, v 
	+ \iO \nabla \ph \cdot \nabla v
	+\iO  \xi v
	+\frac 2 {\thc}\iO  \pi(\ph) v
	- \frac 1{\thc^2} \iO {\dt w}\,\pi(\ph)v = 0\,,
	\\ 
	\label{wf:2}
	& \<\ddt w, v>_{V}
	+\a \iO \nabla (\dt w)\cdot \nabla v
	+ \b \iO \nabla  w  \cdot \nabla v
	+ \iO \pi (\ph) \dt \ph \, v
	=
	\iO \an{u}  v\,,
\end{align}
\pier{are satisfied for every test function $v \in V$ \an{and} almost everywhere in $(0,T)$.}
Moreover, it holds that
\begin{align*}
	\pier{\ph(0)=\ph_0 
	,
	\quad w(0)=w_0 
	,
	\quad \dt w(0)=\an{v_0} 
	.}
\end{align*}
Furthermore, there exists a constant $K_1>0$, \juerg{which depends only} on $\Omega,T,\a,\b,\thc$ and the data of the system, such that
\begin{align}
	\non
	& \norma{\ph}_{\H1 H \cap \L\infty V \cap \L2 {\Hx2}}
	+ \norma{\hat{\gamma}(\ph)}_{\L\infty {\Lx1}}^{1/2}
	\\
	& \quad \label{weaksol:estimate}
	+ \norma{w}_{\H2{\Vp} \cap\W{1,\infty} H \cap \H1 V}
	\leq K_1.
\end{align}
\end{theorem}
Let us emphasize that the above result is very general and includes all of the choices for the potentials introduced in \eqref{Freg}--\eqref{F2ob}.
Besides, notice that the second condition in \eqref{ass:exweak:initialdata} follows from the first one in the case of \eqref{Freg}. In fact,
\pier{we have that} $\hat{\gamma}_{\rm reg}(r) = \mathcal{O}(r^4)$ as $\juerg{|r|} \to \infty$, and
in the three-dimensional case it \pier{turns out} that $\ph_0 \in V \subset \Lx6$. \pier{In view of the regularity of the solution, note that the initial conditions \juerg{make} sense at least in $H$, since, in particular, $\ph \in \C0V $ and $w\in \C1H$ by interpolation properties. Moreover, terms like the last integrals on the 
left-hand sides of \eqref{wf:1} and \eqref{wf:2} are well defined thanks to H\"older's inequality, since $\dt w \in \L2{V}$, $\pi(\ph) \in \L\infty{V} $,  $\dt \ph \in \L2H$, and $ V \subset L^{p}(\Omega)$ for $1\leq p\leq 6$.}

\begin{proof}[Proof of Theorem \ref{THM:EX:WEAK}]
\pier{We proceed by formal estimates, referring\an{, e.g.,} to the papers~\cite{CG1,CGS26} for the details on a regularization and Faedo--Galerkin approximation of a similar but abstract system.}

\noindent
{\bf First estimate:}
\pier{Note that \eqref{sys1} or, more precisely, 
\begin{equation} 
\label{sys1-pier}
\displaystyle \dt\ph-\Delta\ph + \xi + \frac 2 {\thc} \pi(\ph) - \frac 1{\thc^2}{\dt w}\,\pi(\ph) = 0 \quad\mbox{in }\,Q\\,
\end{equation}
with $\xi \in \gamma (\ph)$ \an{almost everywhere} in $Q$, and \eqref{sys2} are the equations related to the variational equalities \eqref{wf:1} and \eqref{wf:2}, respectively.
We test \eqref{sys1-pier}} by $\thc^2 \dt \ph$ and \eqref{sys2} by $\dt w$.
Then we add the resulting equalities and to both sides the term $\frac {\thc^2}2 \big(\norma{\ph(t)}^2 - \norma{\ph_0}^2 \big) = \thc^2 \intQt \ph \, \dt\ph$. \pier{Note that there is a cancellation of two terms. Integrating by parts, we obtain} that
\begin{align*}
	& \thc^2 \I2 {\dt\ph} 
	+ \frac{\thc^2}2 \norma{\ph(t)}^2_V
	+ \thc^2 \iO \hat{\gamma}(\ph(t))
	\\
	&\qquad{}+\frac12 \IO2 {\dt w}
	+ \a \I2 {\nabla (\dt w)}
	+ \frac \b 2 \IO2 {\nabla w}
	\\
	& \quad 
	\pier{{}\leq{}}
	\frac{\thc^2}2 \norma{\ph_0}^2_V
	+\thc^2 \iO \hat{\gamma}(\ph_0)
	+ \frac12 \norma{\an{v_0}}^2
	+\frac \b 2 \norma{\nabla w_0}^2
	\\ & \qquad{}
	-  2 \thc\intQt \pi(\ph) \dt \ph 
	+ \intQt \an{u}\, \dt w
	+\thc^2 \intQt \ph\,\dt\ph.
\end{align*}
The first four terms on the \rhs\ are easily bounded due to the assumption\pico{~\eqref{ass:exweak:initialdata}} \an{on the initial data}.
As for the other three terms, we have, using \an{\eqref{ass:regf:weak},} Young's inequality and the \Lip\ continuity of $\pi$, that
\begin{align*}
	 -2 \thc\intQt  \pi(\ph) \dt \ph 
	 +\thc^2 \iO \ph\,\dt\ph & \leq \frac{\thc^2}2 \I2 {\dt\ph} 
	 + \pier{C}\intQt (|\ph|^2+1),
	 \\
	\intQt \an{u}\, \dt w & \leq 
	\pier{\frac 1 2} \intQt|\dt w |^2 + C.
\end{align*}
\pier{Now, we can apply Gronwall's lemma, which finally} entails that
\begin{align}
	& \norma{\ph}_{\H1 H \cap \L\infty V}
	+ \norma{\hat{\gamma}(\ph)}_{\L\infty {\Lx1}}^{1/2}
	+ \norma{w}_{\W{1,\infty} H \cap \H1 V}
	\leq C. \label{pier1}
\end{align}

\noindent
{\bf Second estimate:}
Next, we \an{take} \pier{an} arbitrary function $v \in \L2 V$ in \pier{\eqref{wf:2}, then use
the linear growth of $\pi$, \Holder's inequality, and the continuous inclusion $V \subset \Lx6$, to infer that}
\begin{align*}
	&\Big| \ioT \< \ddt w, v>_{V} \,\juerg{{\rm d}t}\Big| \\
	&\quad \leq C \ioT \Big( 
	\norma{\nabla (\dt w)}\norma{\nabla v}
	+ \norma{\nabla w}\norma{\nabla v}
	+ \norma{\an{u}} \norma{v}
	\Big)\,\juerg{{\rm d}t}
	+ C \intQ (|\ph|+1)|\dt\ph| |v|
	\\ &\quad
	\leq C \ioT \Big( 
	\norma{\nabla (\dt w)}
	+ \norma{\nabla w}
	+ \norma{\an{u}}
	+ (\norma{\ph}_3+1) \norma{\dt\ph}
	\Big) \norma{v}_V\,\juerg{{\rm d}t}
	\\ & \quad
	\leq C \norma{v}_{\L2 V}.
\end{align*}
Thus, it is a standard matter to \pier{conclude} that
\begin{align}
	\norma{\ddt w}_{\L2 {\Vp}} \leq C. \label{pier2}
\end{align}

\noindent
{\bf Third estimate:}
Next, we notice that \eqref{sys1-pier} can be rewritten as the elliptic equation
\begin{align*}
	 -\Delta \ph + \xi \, \pier{=g , \quad \hbox{with} \quad g \an{{}:={}} - \dt\ph - \frac 2 {\thc} \pi(\ph) + \frac 1{\thc^2}{\dt w}\,\pi(\ph)}
\end{align*}
and $\xi \in \gamma(\ph)$ \an{almost everywhere in $Q$. 
Due to the estimate~\eqref{pier1}\pico{,} $g$ is bounded in $\L2 H$: indeed, it turns out} that $\dt w \in \L2 {\Lx4}$ and $\pi(\ph) \in \L\infty {\Lx4}$.
Thus, \pier{formally testing by $ -\Delta \ph $ and using monotonicity to infer that 
$\int_\juerg{Q} \xi (-\Delta \ph) \geq0$, we find that
\begin{align*}
	\norma{\juerg{\Delta \ph}}_{\L2 H} + \norma{\xi}_{\L2 H}\leq C.
\end{align*}
Then, from \eqref{sys1-pier}, the smooth boundary condition \eqref{sys3} for $\ph$,
 and well-known elliptic regularity results (see, e.g., \cite{Bre}), it follows that}
\begin{align}
	\norma{\ph}_{\L2 {\Hx2}}  \leq C. \label{pier3}
\end{align}
\pier{This ends the proof of \an{the estimate} \eqref{weaksol:estimate}, whence Theorem~\ref{THM:EX:WEAK} is completely proved.}
\end{proof}

\begin{theorem}\label{THM:UQ:WEAK}
Suppose that \ref{ass:1:constants}--\ref{ass:3:pi} hold. 
Then there exists a unique weak solution $(\ph, w, \xi)$ \pier{to the system~\eqref{sys1}--\eqref{sys4}} in the sense of Theorem~\ref{THM:EX:WEAK}. 
Moreover, let us denote by $\{(\ph_i, w_i, \xi_i)\}_{i=1,2}$ \pier{a pair} of weak solutions obtained by Theorem \ref{THM:EX:WEAK} and related to the initial data 
$\{\ph_{\pier{{0,i}}},w_{\pier{{0,i}}},\an{v_{0,i}}\}_{i=1,2}$ and heat sources $\{\an{u}_i\}_{i=1,2}$ fulfilling \eqref{ass:exweak:initialdata} and \eqref{ass:regf:weak}, respectively.
Then it holds that 
\begin{align}
	& \non
	\norma{\ph_1- \ph_2}_{\L\infty H \cap \L2 V}
	+ 
	\norma{w_1- w_2}_{\H1 H \cap \L\infty V	}
	\\ & 
	\quad 	
	\an{\leq K_2 \big ( 
	\norma{\ph_{\pier{{0,1}}}-\ph_{\pier{{0,2}}}}
	{{}+ \norma{w_{\pier{{0,1}}}-w_{\pier{{0,2}}}}_{V}} +\norma{v_{0,1}-v_{0,2}}_{} \big )}
	\non \\ & \quad\quad
	\an{{}+ K_2 \norma{1*(\an{u}_1-\an{u}_2)}_{\L2 H}}
	\label{cd:weak}
\end{align}
\juerg{with a positive constant $K_2$ that depends only on} $\Omega,T,\a,\b,\thc$ and the data of the system.
\end{theorem}

\begin{proof}[Proof of Theorem \ref{THM:UQ:WEAK}]
\pier{We aim to prove} the stability estimate \eqref{cd:weak}. 
This will in turn guarantee the uniqueness of \juerg{weak solutions}.
For convenience, let us set
\begin{align}
	\label{not:diff:1}
	&\ph  : = \ph_1 - \ph_2, 
	\quad 
	w : = w_1 - w_2, 
	\quad 
	\xi:= \xi_1- \xi_2,
	\\
	\label{not:diff:2}
	&
 	\pier{\rho_i := \pi(\ph_i) \quad \text{for } \, i=1,2,
	\quad	
	\rho  := \rho_1- \rho_2,}
	\\	
	& \ph_0   : = \ph_{\pier{{0,1}}} - \ph_{\pier{{0,2}}}, 
	\quad 
	w_0  : = w_{\pier{{0,1}}} - w_{\pier{{0,2}}}, 
	\quad 
	\an{v_0}  : = \an{v_{0,1} - v_{{{0,2}}}}, \quad \an{u}  := \an{u}_1 - \an{u}_2.
	\label{not:diff:3}
\end{align}
Using this notation, we take the difference of the weak formulation \eqref{wf:1}--\eqref{wf:2} written for 
$\{(\ph_i, w_i, \xi_i)\}_{i=1,2}$ and 
\pier{$\{\ph_{\pier{{0,i}}},w_{\pier{{0,i}}}, \pico{v}_{\pier{{0,i}}}, \an{u}_i\}_{i=1,2}$}, 
\pier{obtaining} that the differences fulfill
\begin{align}
	\label{wf:cd:1}
	& \iO \dt\ph \, v 
	+ \iO \nabla \ph \cdot \nabla v
	+\iO  \xi v
	+\frac 2 {\thc}\iO  \rho \, v
	- \frac 1{\thc^2} \iO{\dt w}\,\rho_1 v 
	- \frac 1{\thc^2} \iO {\dt w_2}\,\rho\, v = 0
	,
	\\ 
	\label{wf:cd:2}
	& \<\ddt w, v>_{V}
	+\a \iO \nabla (\dt w)\cdot \nabla v
	+ \b \iO \nabla  w  \cdot \nabla v
	+ \iO \dt (\hat \pi(\ph_1)-\hat \pi(\ph_2)) v
	=
	\iO \an{u}  v,
\end{align}
\pier{for all $v \in V$ and almost everywhere in $(0,T)$. 
Note that, thanks to \ref{ass:3:pi}, we could write the terms $\rho_i \, \dt \ph_i $ appearing in \eqref{sys2}  as $\dt \hat \pi(\ph_i) $, $i=1,2$. Of course, also the initial conditions 
\begin{align}
\label{pier4}
	\ph(0)=\ph_0,\quad w(0)=w_0,\quad \dt w(0)=\an{v_0},\quad \mbox{ hold a.e. in }\,\Omega.
\end{align}
First, we add the term $\int_\Omega \ph \, v $ to both sides of \eqref{wf:cd:1}, then take 
$v= \ph$  and integrate with respect to time. We deduce that
\begin{align}
	& \frac {1}2 \IO2 {\ph}
	+ \int_0^t \norma{\ph(s)}_V^2 \, \ds
	+ \intQt \xi\,\ph \non
	\\ & \quad 
	=
	\frac {1}2 \norma{\ph_0}^2
	+ \intQt  \Big( \ph  - \frac{2}{\thc} \rho \Big) \ph
	+ \frac {1}{\thc^2}\intQt {\dt w}\,\rho_1 \, \ph
	+ \frac {1}{\thc^2}\intQt {\dt w_2}\,\rho \,\ph
	\label{pier3-1}
\end{align}
for all $t\in [0,T]$.
Due to the monotonicity of $\gamma$, we immediately conclude that the third term on the \lhs\ is nonnegative.
Using the \Lip\ continuity of $\pi$ along with the regularities $\dt w_i \in \L\infty H \cap \L2 V$, $\ph_i \in \H1 H \cap \L\infty V$, $i=1,2$,  we infer from Theorem~\ref{THM:EX:WEAK} that
\begin{align*}
&{} \intQt  \Big( \ph  - \frac{2}{\thc} \rho \Big) \ph
 \leq C \I2 \ph\,,
 \end{align*} 
 and, with the help of H\"older's inequality and of the continuous embedding $ V \subset L^4(\Omega)$,
\begin{align*} 
&\frac {1}{\thc^2}\intQt {\dt w}\,\rho_1 \ph
 \leq  C\int_0^t \norma{\dt w} \Big( \norma{\ph_1}_4 +1\Big) \norma{\ph }_4\,\juerg{{\rm d}s}
\\ 
&\quad
\leq  C \Big( \norma{\ph_1}_{\L\infty V} +1\Big)\int_0^t \norma{\dt w} \, \norma{\ph }_V \,\juerg{{\rm d}s}
  \leq  \frac 1 4 \int_0^t\norma{\ph }_V^2 \,\juerg{{\rm d}s}+ D_1 \int_{Q_t}  |\dt w|^2 , 
\end{align*}  
where $D_1$ is a computable and by now fixed constant.
Moreover, we have that 
\begin{align*}  
&\frac {1}{\thc^2}\intQt {\dt w_2}\,\rho \,\ph	 
\leq C\int_0^t \norma{\dt w_2}_4 \, \norma{\ph} \norma{\ph }_4 \,\juerg{{\rm d}s}
	\\  
&\quad\leq	C \int_0^t \norma{\dt w_2}_V \, \norma{\ph } \, \norma{\ph }_V	\,\juerg{{\rm d}s}
\leq  \frac 1 4 \int_0^t\norma{\ph }_V^2 \,\juerg{{\rm d}s}+C  \int_0^t \norma{\dt w_2}_V^2 \, \norma{\ph }^2 
\,\juerg{{\rm d}s},	
\end{align*}
where the function $ t\mapsto \norma{\dt w_2(t) }_V^2$ belongs to $L^1(0,T)$ \an{due to Theorem~\ref{THM:EX:WEAK}}. Therefore, collecting the above estimates,  it follows from \eqref{pier3-1} that 
\begin{align}
	& \frac {1}2 \IO2 {\ph}
	+ \frac {1}2\int_0^t \norma{\ph(s)}_V^2 \, \ds
	\non \\ & \quad 
	\leq
	\frac {1}2 \norma{\ph_0}^2
+  C  \int_0^t \Big(1 +\norma{\dt w_2}_V^2 \Big)\, \norma{\ph }^2\,\juerg{\ds} 
+ D_1 \int_{Q_t}  |\dt w|^2	.
	\label{pier3-2}
\end{align}}%
Next, we integrate  \eqref{wf:cd:2} with respect to time using \eqref{pier4}, then take $v= \dt w$, and integrate once more over $(0,t)$\an{, for an arbitrary $t \in [0,T]$}. \juerg{Addition of the terms 
$\frac {\a}2 \big(\norma{w(t)}^2 - \norma{w_0}^2 \big) = \a \intQt w \, \dt w$ to both sides} leads to
\begin{align}
	& \I2 {\dt w }
	+ \frac\a2 \norma{w(t)}^2_V
	=
	\intQt {\an{v_0}}\, \dt w
	+ \intQt (\hat \pi(\ph_{\pier{{0,1}}})-\hat \pi(\ph_{\pier{{0,2}}}))\dt w
	\non \\ & \quad 
	+ \a\intQt { \nabla w_0} \cdot \nabla(\dt w)
	+\frac\a2 \norma{w_0}
	\an{{}-{}} \b \intQt (1 * \nabla w) \cdot \nabla(\dt  w)
	\non \\ & \quad 
	- \intQt (\hat \pi(\ph_1)-\hat \pi(\ph_2)) \dt w
	+ \intQt (1*\an{u} ) \dt w
	+\a \intQt w \, \dt w.
	\label{pier3-3}
\end{align}
We estimate each term on the \rhs\ individually. 
Let us recall that the mean value theorem and the Lipschitz continuity of $\pi$ yield 
\an{the existence of a positive constant $C$ such}
that
\begin{align}\label{pier5}
	\left| \hat \pi(r)-\hat \pi(s) \right| 
	\leq C (|r| + |s| + 1) |r-s|\pico{\quad \hbox{for all } \, r,\, s \in \erre.}
\end{align}  
By Young's inequality, we easily have 
\begin{align*}
	& \intQt {\an{v_0}}\,\dt w
	\leq 
	\frac18 \I2 {\dt w}
	+ C \norma{{\an{v_0}}}^2.
\end{align*} 
Using integration over time, \Holder's inequality, \eqref{pier5}, and the continuous embedding $V \subset \Lx4$, we find that
\begin{align*}
	& \intQt (\hat \pi(\ph_{\pier{{0,1}}})-\hat \pi(\ph_{\pier{{0,2}}}))\dt w
	= \iO (\hat \pi(\ph_{\pier{{0,1}}})-\hat \pi(\ph_{\pier{{0,2}}})) (w(t) - w_0) 	
	\\ & \quad \leq 
	C \big\| |\ph_{\pier{{0,1}}}| + |\ph_{\pier{{0,2}}}| + 1\big\|_4 \,
	\norma{\ph_{\pier{{0,1}}}- \ph_{\pier{{0,2}}}} 
	\big(\norma{w(t)}_4 + \norma{ w_0}_4	\big)
	\\ & \quad \leq	
	C\big(\norma{\ph_{\pier{{0,1}}}}_V+\norma{\ph_{\pier{{0,2}}}}_V+1\big) 	
	\norma{ \ph_0}\big(\norma{w(t)}_V + \norma{ w_0}_V	\big)	
	\\ & \quad \leq		
    \frac\a 8 \big(\norma{w(t)}^2_V + \norma{ w_0}^2_V	\big) + C \big(\norma{\ph_{\pier{{0,1}}}}_V^2+\norma{\ph_{\pier{{0,2}}}}_V^2 +1\big) \norma{ \ph_0}^2  .     
\end{align*}
Next, the third term on the \rhs\ of \eqref{pier3-3} can be bounded as
\begin{align*}
		\a \intQt {\nabla w_0} \cdot \nabla(\dt w)
		= \a \iO \nabla w_0 \cdot (\nabla w(t) - \nabla w_0)
		\leq 
		\frac\a 8  \IO2 {\nabla w}
		+ C \norma{\nabla w_0}^2.
\end{align*}
Then, by using the identity
\begin{align*}
	\intQt (1*\nabla w) \cdot \nabla (\dt w )
	= \iO (1*\nabla w(t)) \cdot \nabla w(t)
	- \intQt  |\nabla w|^2, 
\end{align*}
the fact that \,$\Vert 1*\nabla w(t) \Vert^2 \leq \Big( \int_0^t \Vert \nabla w \Vert\Big)^2\leq T\intQt |\nabla w|^2$, and Young's inequality, we infer that
\begin{align*}
	- \b \intQt (1 * \nabla w) \cdot \nabla (\dt w )
	\leq \frac\a 8 \IO2 {\nabla w}
	+ C \I2 {\nabla w}.
\end{align*}
To handle the sixth term on the \rhs\ \an{of \eqref{pier3-3}}, we owe once more to \eqref{pier5} and 
the continuous and compact embedding $V \subset \Lx{p}$, $1\leq p <6 $.
By the H\"older and Young inequalities, and thanks to \eqref{weaksol:estimate} and 
the Ehrling lemma (see, e.g., \cite[Lemme~5.1, p.~58]{Lions}), we can deduce that
\begin{align*}
	&- \intQt (\hat \pi(\ph_1)-\hat \pi(\ph_2)) \dt w
	\leq 
	C \iot \big\| |\ph_{\pier{{1}}}| + |\ph_{\pier{{2}}}| + 1 \big\|_4 \,
	\norma{\ph_{\pier{{1}}}- \ph_{\pier{{2}}}}_4 \norma{\dt w} \juerg{\,\ds}
	\\ &  \quad 
	\leq 
	\frac 18 \I2 {\dt w }
	+ {C}(\norma{\ph_1}_{\L\infty{V}}^2+\norma{\ph_2}_{\L\infty{V}}^2+1) 
	\int_0^t  \|\ph\|_4^2 \juerg{\,\ds}
	\\ &  \quad 
	\leq 
	\frac 18 \I2 {\dt w } + \delta \int_0^t  \|\ph\|_V^2	\juerg{\,\ds}+ C_\delta \I2 \ph	\,,
\end{align*}
for \juerg{any} positive coefficient $\delta$ \juerg{(yet to be chosen)}. 
Lastly, Young's inequality easily produces
\begin{align*}
	\intQt (1*\an{u} ) \dt w
	+	\a \intQt w \, \dt w 
	& 
	 \leq \frac 14  \I2 {\dt w}
	 + C\I2 {1*\an{u}}
	+ C \I2 w.
\end{align*}
Thus, in view of \eqref{pier3-3}, upon collecting the above computations, we realize that 
\begin{align}
	& \an{\frac 12} \I2 {\dt w }
	+ \frac\a 8 \norma{w(t)}^2_V
		\non \\ & \quad 
		\leq 
	C \norma{{\an{v_0}}}^2
	+ C\big(\norma{\ph_{\pier{{0,1}}}}^2_V+\norma{\ph_{\pier{{0,2}}}}^2_V+1\big) \norma{ \ph_0}^2	+ C \norma{w_0}^2_V
	\non \\ & \quad\quad 
	+ \delta \int_0^t  \|\ph\|_V^2 \juerg{\,\ds}	+ C_\delta \I2 \ph	
	 + C\I2 {1*\an{u}}
	+ C \int_0^t \norma{w}_V^2 \juerg{\,\ds}.	
	\label{pier3-4}
\end{align}

At this point, we multiply \eqref{pier3-4} by $\an{4}D_1 $ and add it to \eqref{pier3-2};
then, fixing $\delta>0 $ such that $\an{4}D_1 \delta <1/2$, and applying the Gronwall lemma, \juerg{we obtain} the estimate 
\begin{align*}
	& \norma{\ph}_{\L\infty H \cap \L2 V}
	 + \norma{w}_{\H1 H \cap \L\infty V}
	 \\ & \quad
	 \leq C ({\norma{ \ph_0}}  +\norma{w_0}_{V} + \norma{\an{v_0}} + \norma{1*\an{u}}_{\L2 H}),
\end{align*}
where $C$ depends also on $\norma{\ph_{\pier{{0,i}}}}_V$, $i=1,2$. Due to our notation in 
\eqref{not:diff:1}--\eqref{not:diff:3}, this is actually \eqref{cd:weak}, and the proof \an{of Theorem~\ref{THM:UQ:WEAK}} is complete.
\end{proof}

To improve the regularity results of Theorem \ref{THM:EX:WEAK}, as well as the stability estimate \eqref{cd:weak}, we are forced to require more regularity \an{on structural elements}, \juerg{in particular, for} the nonlinearity $\hat \gamma$. In the following lines, we 
\juerg{state} general conditions under which we are able to extend the existence and uniqueness results to a stronger framework.
\begin{enumerate}[label={\bf B\arabic{*}}, ref={\bf B\arabic{*}}]
\item \label{ass:strong:gamma:1}
	There exists an interval $\,(r_-,r_+)\,$ with $\,-\infty\le r_-<0<r_+\le +\infty\,$ such that
the restriction of $\hat \gamma$ to $\,(r_-,r_+)\,$ belongs to $\,{C^2}(r_-,r_+)$.
\pier{Thus, $\gamma$ coincides with the derivative of $\hat \gamma$ in $(r_-,r_+)$.}

\item \label{ass:strong:gamma:2}
	It holds that\,\, \juerg{$\lim_{r\searrow r_-} \gamma(r)=-\infty$ \,and\, $\lim_{r\nearrow r_+}\gamma(r)=+\infty$}.
\item \label{ass:4:pidiff}
	\pier{$\gamma \in {C^2}(r_-,r_+)$	and $ \pi \in C^2 (\erre)$.}
	\end{enumerate}
Notice that \ref{ass:strong:gamma:1}--\ref{ass:4:pidiff} are fulfilled by the regular and the logarithmic potentials \eqref{Freg} and \eqref{Flog}, whereas the double obstacle nonlinearity \eqref{F2ob} is no longer allowed. Again, we remark that, due to \ref{ass:strong:gamma:1}, we no longer need to consider any selection $\xi \in \partial \gamma(\ph)$ as $\gamma  = \hat \gamma \,'$ \pier{in $\,(r_-,r_+)\,$}. This also entails that \eqref{sys1} 
\juerg{becomes} an equality. 

\pier{The next result dealing with regularity of the solution does not need
the condition~\ref{ass:4:pidiff}.}
\begin{theorem}
\label{THM:EX:STRONG}
\pier{Assume \an{that} \ref{ass:1:constants}--\ref{ass:3:pi} and \ref{ass:strong:gamma:1}--\ref{ass:strong:gamma:2}} \an{are fulfilled}.
Furthermore, let the heat source $\an{u}$ \pier{fulfill} \eqref{ass:regf:weak}, and \pier{let} the initial data, in addition to \eqref{ass:exweak:initialdata}, satisfy
\begin{align}\label{ass:exstrong:initialdata}
	\pier{{}\ph_0 \in W,\quad   \an{v_0} \in V, \quad 
	\ph_0' := \Delta \ph_0 - \gamma(\ph_0) - \tfrac 2{\thc} \pi(\ph_0) + \tfrac 1 {\thc^2}\an{v_0} \pi(\ph_0)   \in H. {}}
\end{align}
Then there exists a strong solution $(\ph, w)$ to system \Sys\ in the sense that 
\begin{align}
	\label{reg:strong:ph}
	\ph &\in {\W{1,\infty} H \cap \H1 V} \cap \L\infty{W},
	\\ 
	\label{reg:strong:w}
	w & \in {\H2 H \cap \W{1,\infty} V \cap \H1 {W}},
\end{align}
and that \an{the equations} \Sys\ \pier{are} fulfilled almost everywhere in $Q$, on $\Sigma$, or in $\Omega$, respectively.
In addition, assume that the heat source $\an{u}$ fulfills 
\begin{align}
	\label{ass:regf:strong}
	\an{u} \in \L\infty H 
\end{align}
and that 
\begin{align}
	\pier{w_0 , \, \an{v_0} \in \Lx\infty},
	\quad 
	r_-<\min_{x\in\overline{\Omega}}\,\vp_0(x) \le
	\max_{x\in\overline{\Omega}}\,\vp_0(x)<r_+.
	\label{strong:sep:initialdata}
\end{align}
Then it holds that
\begin{align*}
	\dt w \in L^\infty(Q)\,,
\end{align*}
and the phase variable $\ph$ enjoys the so-called {separation property}, \juerg{which means that}
there exist two values $r_*,r^*$, depending \juerg{only} 
on $\Omega,T,\a,\b, \thc$ and the data of the system, such that
\begin{align}
	\label{separation}
	r_-  <r_* \leq {\ph} \leq  r^* < r_+ \quad \hbox{a.e. in } Q.
\end{align}
Furthermore, there exists a constant $K_3>0$ such that
\begin{align}\label{reg:strong}
	\non
	& \norma{\ph}_{{\W{1,\infty} H \cap \H1 V} \cap \L\infty{\Hx2}}
	\\ & \quad 
	+ \norma{w}_{\H2 H \cap \W{1,\infty} V \cap \H1 {\Hx2}}
	+ \norma{\dt w}_{L^\infty(Q)}
	\leq K_3.
\end{align}
\end{theorem}
Here, \pier{we point out that the regularities in \eqref{reg:strong:w} imply $w \in C^0({\ov Q})$
thanks to the Sobolev embedding results.} \juerg{Moreover, since the embedding $W\subset C^0(\ov \Omega)$
is compact, it follows from \cite[Sect.~8, Cor.~4]{Simon} that also $\ph\in C^0(\ov Q)$. In
particular, the separation property \eqref{separation} is valid even pointwise in $\ov Q$.} 

\begin{proof}[Proof of Theorem \ref{THM:EX:STRONG}]
\pier{
\an{In what follows, we}
perform the estimate directly on \an{the system} \Sys\ underlying that now also 
equation~\eqref{sys1} \an{turns to} an equality as 
$\gamma (\cdot) = \hat \gamma^{\,\prime} (\cdot)$ is single valued. \an{A} rigorous proof would need \an{some} approximation, but please take into account that we already have proved the existence and uniqueness of the weak solution.}

\noindent
{\bf First estimate:}
To begin with, we \pier{formally} differentiate \eqref{sys1} with respect to 
time and \pier{multiply} the resulting identity by $\thc^2 \dt \ph$\pier{;
then we add \eqref{sys2} tested by $\ddt w$, and integrate over $Q_t$. Note that a 
cancellation occurs and that, after some rearrangements, one obtains}
\begin{align*}
	& \frac{\thc^2}2 \IO2 {\dt\ph} 
	+ {\thc^2} \I2 {\nabla (\dt \ph)}	
	+ \thc^2 \intQt {\gamma}'(\ph)|\dt\ph|^2
	+ \I2 {\ddt w}
	+ \frac \a 2 \IO2 {\nabla (\dt w)}
	\\
	& \quad 
	\pier{{}\leq{}}
	\frac{\thc^2}2 \norma{\ph_0'}^2
	+\frac \a 2 \norma{\nabla \an{v_0}}^2
	-  2 \thc \intQt \pi'(\ph)|\dt\ph|^2
	+ \intQt \dt w \, \pi'(\ph) |\dt \ph|^2	
	\\ & \qquad 
	- \b \intQt \nabla w \cdot \nabla (\ddt w)
	+ \intQt \an{u} \,\ddt w.
\end{align*}
\pier{Owing to 
the monotonicity of $\gamma$, we infer that the third term on the \lhs\ is nonnegative.}
The first two terms on the \rhs\ are controlled due to \an{the conditions \pico{\eqref{ass:exstrong:initialdata}} on the initial data}.
\pier{As for the third term on the \rhs, we note that $\dt \pi (\ph) $
 makes sense as $ \pi'(\ph)\dt \ph $, \juerg{in view of the global} Lipschitz continuity of $\pi$. Now,  
we use the boundedness of $\pi'$ and
estimate~\eqref{weaksol:estimate}, obtaining} \an{that}
\begin{align*}
	-  2 \thc \intQt \pi'(\ph)|\dt\ph|^2
	\leq C \I2 {\dt \ph} \pier{{}\leq C}. 
\end{align*}
\pier{Next, as $V\subset \Lx4 $ with compact embedding, we employ H\"older's inequality, \eqref{weaksol:estimate}, and Ehrling's lemma, to deduce} that
\begin{align*}
	\intQt \dt w \, \pi'(\ph) |\dt \ph|^2	 &\leq C \iot \norma{\dt w} \norma{\dt \ph}_4^2 \juerg{\,\ds}
	\\
	&
	\leq C \iot \norma{\dt \ph}_4^2 \juerg{\,\ds}
	\leq \frac{\thc^2}2 \I2 {\nabla (\dt \ph)} + C \I2 {\dt \ph}.
\end{align*}
The fifth term on the \rhs\ can be controlled by integrating \pier{by parts} and using the above estimate along with Young's inequality and assumptions \eqref{ass:exstrong:initialdata}, \pier{so that}
\begin{align*}
	&- \b \intQt \nabla w \cdot \nabla (\ddt w) 
	 \\ & \quad 
	 = 
	\b \I2 {\nabla(\dt w)}
	- \b \iO \nabla w(t) \cdot \nabla (\dt w(t))
	+ \b \iO \nabla w_0 \cdot \nabla \an{v_0}
	\\ &\quad  \leq
	\b \I2 {\nabla(\dt w)}
	+ \frac \a4 \IO2 {\nabla( \dt w)} 
 	+ C \norma{w}_{\L\infty V}^2 
 	+ C (\norma{w_0}_V^2 + \norma{\an{v_0}}_V^2).
\end{align*}
\pier{Finally, the last term can be easily handled by Young's inequality, namely,}
\begin{align*}
	\intQt \an{u}\, \ddt w \leq \frac 12 \I2 {\ddt w} + \pier{{} \frac12 \norma{\an{u}}_{\L2H}^2 .}
\end{align*}
Hence, upon collecting the above computations, the Gronwall lemma yields that 
\begin{align}
	& \norma{\ph}_{\W{1,\infty} H \cap \H1 V}
	+ \an{\norma{w}_{\H2 H \cap \W{1,\infty} V}}
	\leq C.
	\label{pier6}
\end{align}

\noindent
{\bf Second estimate:}
\pier{By comparison} in equation \eqref{sys1}, we deduce that 
\begin{align*}
	\an{\norma{-\Delta \ph + \gamma(\ph)}}_{\L\infty H} \leq C . 
\end{align*}
\pier{Then, arguing as in the proof of Theorem~\ref{THM:EX:WEAK} (cf.~the Third estimate there),
and using the elliptic regularity theory, we infer that}
\begin{align}
	\label{est:phinftyH2}
	\norma{\ph}_{\L\infty {\Hx2}} + \norma{\gamma(\ph)}_{\L\infty {H}}\leq C .
\end{align}

\noindent
{\bf Third estimate:}
We then rewrite \eqref{sys2} as a parabolic \pier{equation in the new variable $y := \a \dt w + \b w$. Thanks to \an{equations} \eqref{sys3}--\eqref{sys4}, we have that}
\begin{align}\label{parab:sys}
  \begin{cases}
 	\displaystyle \an{\tfrac 1\a} \dt y -\Delta y = g : = \an{u} - \pi (\ph) \dt \ph + \an{\tfrac \b\a} \dt w \quad & \text{in}\,\, Q,
 	\\[2mm]
 	\dn y = 0 \quad & \text{on}\,\, \Sigma,
 	\\[1mm]
 	y(0) = y_0 :=\a \hskip1pt\an{v_0} + \b w_0 \quad & \text{in}\,\, \Omega.
  \end{cases}
\end{align}
By analyzing system \eqref{parab:sys}, we realize that $g \in \L2 H$ and $y_0 \in V$, 
so that \pier{the parabolic regularity theory} entails that
\begin{align}
	\norma{y}_{\H1 H \pier{{}\cap \L\infty V{}}\cap \L2 {\Hx2}} \leq C.
\label{pier6-1}
\end{align}
In fact, \pier{since the ODE relation $\a \dt w + \b w = y $ holds true in $ Q,$ then 
\begin{equation}
w(t) = e^{-\beta t /\alpha} w_0 + \frac 1 \a \int_0^t e^{-\beta (t -s) /\alpha}y(s) \pico{\ds} , \quad t\in [0,T]. \label{pier6-9}
\end{equation}
Thus, $w$ and its derivative $\partial_t w$ possess the same regularity as $y$ and satisfy estimates like \eqref{pier6-1}, where the constant on the \rhs\ has the same dependencies. 
Therefore, we} eventually conclude that
\begin{align}
	\norma{w}_{\H2 {H} \pier{{}\cap \W{1,\infty} V{}}\cap \H1 {\Hx2}} \leq C.
		\label{pier7}
\end{align}

\noindent 
{\bf Fourth estimate:}
Let us consider again system \eqref{parab:sys}. Due to the above estimates and to \eqref{ass:regf:strong}\pier{, we have that $g $ is bounded in} $\L\infty H$. \pier{Thanks
to \eqref{ass:exweak:initialdata}, \eqref{ass:exstrong:initialdata}, and the first condition in \eqref{strong:sep:initialdata}, it turns out that the initial datum 
$y_0$ is bounded in $ V \cap \Lx\infty$.
Hence, \juerg{an} application of \cite[Thm.~7.1, p.~181]{LSU} yields that
\begin{align*}
	\norma{y}_{L^\infty(Q)} 
	 = \norma {\a \dt w + \b w }_{L^\infty(Q)} \leq C.
	 \end{align*}
\juerg{Moreover}, arguing as above, this in particular leads to  
\begin{align}
	\label{reg:dtwinfty}
	\norma{w}_{L^\infty(Q)} + \norma { \dt w }_{L^\infty(Q)} \leq C. 
\end{align}
As a consequence, by virtue of \eqref{pier6}, \eqref{est:phinftyH2},  and \eqref{pier7},
 the estimate \eqref{reg:strong} eventually follows.}

\noindent
{\bf Separation property:}
Now, with the help of the regularity result proved above, we are in a position to prove the 
separation property for the phase variable $\ph$. 
This can be \pier{shown} by following the same lines of \pier{argumentation} \juerg{as}
in \cite[Proof of Theorem 2.2]{CSS1} \pier{(see also \cite{CSS1-corr})}. \pier{Observe that $\ph$ is bounded in $L^\infty (Q)$
due to \eqref{est:phinftyH2} and the Sobolev embedding $\Hx2 \subset \Lx\infty$
\juerg{(as noted above, we even have $\ph\in C^0(\ov Q)$)}. Hence,
if we} rewrite \eqref{sys1} as
\begin{align}
 	\dt\ph-\Delta\ph + \gamma(\ph) = g 
	\pier{, \quad \hbox{where now} \quad 
	 g \an{:}= - \an{\tfrac 2 {\thc} \pi(\ph) +\tfrac 1{\thc^2}}{\dt w}\,\pi(\ph),}
	\label{pier8}
\end{align}
\pier{then it turns out that $g$ is bounded in $L^\infty(Q)$, due to \ref{ass:3:pi} and \eqref{reg:dtwinfty}.}
This entails the existence of a positive constant $g^*$ for which $\norma{g}_{L^\infty(Q)}\leq g^*$.
Furthermore, the growth assumptions \ref{ass:strong:gamma:1}--\ref{ass:strong:gamma:2} 
ensure the existence of some constants $r_{*}$ and $r^{*}$ such that $ r_- < r_{*}\leq r^{*} < r^+ $
and
\begin{gather}
	\label{sep:cond:1}
	 r_{*} {{}\leq{}} \juerg{\min_{x\in \ov\Omega}}\,\, \ph_0(x), \quad r^{*} {{}\geq{}} 
	\juerg{\max_{x\in \ov\Omega}}\,\, \ph_0(x),
	\\
	\gamma(r)
	 + g^*\leq 0 
	\quad \forall r \in (r_-,r_{*}), \quad 
	\gamma(r
	)
	 - g^*\geq 0 
	\quad \forall r \in (r^{*},r_+).
	\label{sep:cond:2}
\end{gather}
Then, \pier{if} we set $\lambda= (\ph-r^*)^+$, where $(\hskip1pt\cdot\hskip1pt)^+ := \max \{\hskip1pt\cdot\hskip1pt, 0\}$ denotes the positive part function, \pier{and multiply equation
\eqref{pier8} by $\lambda$, \juerg{then integration} over $Q_t$ and by parts leads to}
\begin{align*}
	\pcol{{}\frac 12{}} \norma{\lambda(t)}^2 + \I2 {\nabla \lambda}
	+ \intQt (\gamma(\ph)-g)\lambda
	= 0,
\end{align*}
\pier{for all $t\in [0,T]$, 
where we also applied \eqref{sep:cond:1}} to conclude that $\lambda(0)=0$.
Moreover, \eqref{sep:cond:2} yields that the last term on the \pcol{\lhs}\ of the above identity is nonnegative, so that it follows $\lambda = (\ph-r^*)^+ =0$,
which means that $\ph \leq r^*$ almost everywhere {in} $Q$.
The same argument can be applied with the choice $\lambda = -(\ph-r_*)^- $, with $(\hskip1pt\cdot\hskip1pt)^-:= - \min \{0,\hskip1pt\cdot\hskip1pt\}$, to derive the other bound $\ph \geq r_*$ almost everywhere {in} $Q$.
Thus, we end up with \pier{the property~\eqref{separation}
and conclude the proof.}
\end{proof}

Finally, \pier{in the more regular framework we can provide a refined continuous dependence result that complements Theorem~\ref{THM:UQ:WEAK}.}
\begin{theorem}\label{THM:UQ:STRONG}
Suppose that \ref{ass:1:constants}--\ref{ass:3:pi} and \ref{ass:strong:gamma:1}--\ref{ass:4:pidiff} hold. 
\pier{Denote by $\{(\ph_i, w_i)\}_{i=1,2}$ two pairs of strong solutions obtained by Theorem~\ref{THM:EX:STRONG} in correspondence with the initial data $\{\ph_{\pier{{0,i}}},w_{\pier{{0,i}}},\an{v_{0,i}}\}_{i=1,2}$ fulfilling \eqref{ass:exweak:initialdata}, \eqref{ass:exstrong:initialdata}, \eqref{strong:sep:initialdata},
and heat sources $\{\an{u}_i\}_{i=1,2}$ as in \eqref{ass:regf:strong}.}
Then it holds that 
\begin{align}\label{cd:strong}
	\non
	&\pier{ \norma{\ph_1- \ph_2}_{\W{1,\infty} H \cap \H1 V \cap \L\infty {W}}} 
	+ 
	\non \norma{w_1- w_2}_{ \H2 H \cap \W{1,\infty} V \cap \H1 {\pier{W}}}
	\\ & \quad  \non
	\leq K_4 \big ( 
	\norma{\ph_{\pier{{0,1}}}-\ph_{\pier{{0,2}}}}_{\pier W}	
	\pier{{}+ \norma{w_{\pier{{0,1}}}-w_{\pier{{0,2}}}}_{V}	
		+ \norma{\an{v_{0,1}-v_{0,2}}}_{V}{}}	
	\big)
	\\ & \qquad
	\an{+ K_4 \norma{\an{u}_1-\an{u}_2}_{\L2 H}	,}
\end{align}
with a positive constant $K_4$ \juerg{ that depends only} on $\Omega,T,\a,\b,\thc$ and  the data of the system.
\end{theorem}

\begin{proof}[Proof of Theorem \ref{THM:UQ:STRONG}]
First, let us recall the notation introduced in \eqref{not:diff:1}--\eqref{not:diff:3} and again consider the variational system \eqref{wf:cd:1}--\eqref{wf:cd:2}. 
Now, owing to the regularity assumption \ref{ass:strong:gamma:1}, we have \pier{$\xi_i= \gamma(\ph_i)$ for $i=1,2$.} 
Moreover, the separation property \eqref{separation} enjoyed by \pier{both $\ph_i$, $i=1,2$,} 
combined with \ref{ass:strong:gamma:1}--\ref{ass:strong:gamma:2}, yields that $\gamma$ is \Lip\ continuous 
\pier{when restricted to $[r_*,r^*]$.} Besides, due to the improved regularity at disposal, we \juerg{may
now express} the difference \pier{$\dt (\hat \pi(\ph_1)-\hat \pi(\ph_2))$ in \eqref{wf:cd:2} as $\rho_1 \dt \ph_1 -\rho_2 \dt \ph_2 =\rho \dt \ph_1 + \rho_2 \dt \ph$.}

Let us now move on checking \pico{the estimate} \eqref{cd:strong}.

\noindent
{\bf First estimate:}
We test \eqref{wf:cd:1} by $\dt \ph$, \eqref{wf:cd:2} by $\dt w $\juerg{, add the resulting identities,} and integrate 
\pier{over $(0,t)$} to infer that
\begin{align}
	& \I2 {\dt \ph}
	+\frac 12 \IO2 {\nabla \ph}
	+ \frac 12 \IO2 {\dt w}
	+ \a \I2 {\nabla (\dt w)}	
	+ \frac \b2 \IO2 {\nabla w}	
	\non \\ & \quad =
	\frac 12 \norma{\nabla \ph_0}^2
	+ \frac 12 \norma{\an{v_0}}^2
	+ \frac \b2 \norma{\nabla w_0}^2
	- \intQt \big(\gamma(\ph_1)-\gamma(\ph_2)\big) \dt \ph
	\non \\ & \qquad 	
	{}- \frac 2 {\thc} \intQt \pier{\rho}\, \dt \ph
	+ \frac 1{\thc^2} \intQt \dt w \, \pier{\rho_1} \, \dt\ph
	+ \frac 1{\thc^2} \intQt  \dt w_2 \, \pier{\rho} \, \dt\ph
	\non\\ & \qquad 	
	{}- \intQt \pier{\rho} \, \dt \ph_1 \, \dt w
	- \intQt \pier{\rho_2} \, \dt \ph \, \dt w
	+ \intQt \an{u} \, \dt w.
	\label{pier9}
\end{align}
The fourth, fifth, and last \juerg{terms} on the \rhs\ can be easily handled using Young's inequality and 
the \Lip\ continuity of $\pi$ and $\gamma$, namely,
\begin{align*}
	& - \intQt \big(\gamma(\ph_1)-\gamma(\ph_2)\big) \dt \ph
	- \frac 2 {\thc} \intQt \pier{\rho}\, \dt \ph
	+ \intQt \an{u} \,\dt w
	\\ & \quad 
	\leq 
	 \pier{{} \frac 1 4 \I2 {\dt \ph} {}}
	+ C \intQt (|\ph|^2 + |\an{u}|^2 + |\dt w|^2). 
\end{align*}
Due to Theorem~\ref{THM:EX:STRONG}, we have that $\ph_i$\pier{, and consequently $\rho_i$,} are uniformly bounded in $L^\infty(Q)$ for $i=1,2$, so that also the sixth and ninth terms can be easily controlled in a similar fashion as
\begin{align*}
	&\frac 1{\thc^2} \intQt\dt w \, \pier{\rho_1} \, \dt\ph	
	- \intQt \pier{\rho_2} \, \dt \ph \, \dt w
	\\
	&\quad{}\pier{{}\leq \frac 1 4  \I2 {\dt \ph} 
	+ C \big(\norma{\rho_1}_{L^\infty(Q)}^2 +\norma{\rho_2}_{L^\infty(Q)}^2\big)} \I2 {\dt w}.
\end{align*}
As for the remaining two terms, we recall that \pier{$\norma{\dt \ph_i}_{\L\infty H}$ and 
$\norma{\dt w_i}_{\L\infty V}$  are bounded for $i=1,2$, so that the H\"older and Young 
inequalities and the continuous embedding $V \subset \Lx{4}$ imply that
\begin{align*}
	& \frac 1{\thc^2} \intQt  \dt w_2 \, \pier{\rho} \, \dt\ph
	- \intQt \pier{\rho} \, \dt \ph_1 \, \dt w	
	\\ 
	& \quad \leq
	C \int_0^t  \norma{\dt w_2}_4 \, \norma{\ph}_4 \, \norma{\dt\ph} \juerg{\,\ds}
	+ C \int_0^t \, \norma{\ph}_4 \, \norma{\dt\ph_1} \, \norma{\dt w}_4	\juerg{\,\ds}
	\\ 
	& \quad \leq
	\pier{\frac 1 4}  \I2 {\dt \ph} 
	+ C\,\norma{\dt w_2}_{\L\infty V}^2	
	\int_0^t \norma{\ph}_V^2 \juerg{\,\ds}
	\\ & \qquad  
	+ \frac \alpha 2  \intQt \big( |\dt w|^2 + |\nabla (\dt w)|^2 \big)
	+ C\, \norma{\dt \ph_1}_{\L\infty H}^2	
	\int_0^t \norma{\ph}_V^2 \juerg{\,\ds}. 	
\end{align*}}%
At this point, we can collect the above estimates and combine them with \eqref{pier9}.
Then we either apply the Gronwall lemma or take advantage of the \pier{already shown
inequality~\eqref{cd:weak}} to bound the \rhs. Thus, we arrive at
\begin{align}\label{pier10}
	\non
	& \norma{\ph}_{\H1 H \cap \L\infty V}
	+ 
	\norma{w}_{\W{1,\infty} H \cap \H1 {V}}
	\\ & \quad 
	\leq C \big ( 
	\norma{\ph_0}_{V}	
	+ \norma{w_{0}}_{V}	
		+ \norma{\pico{v}_{0}}_{}	
	+ \norma{\an{u}}_{\L2 H}	
	\big).
\end{align}

\noindent 
{\bf Second estimate:}
\pier{Arguing as in \eqref{parab:sys},
 we can rewrite \eqref{wf:cd:2} as a parabolic system in the variable 
$y= \a \dt w + \b w$ with source term $\an{g} := \an{u} - \rho \dt \ph_1 - \rho_2 \dt \ph + \tfrac \b\a \dt w$. 
Since 
$$  \norma{\rho \,\dt \ph_1}^2_{\L2H} \leq C \int_0^T \norma{\ph }_4^2 \,  \norma{\dt \ph_1}_4^2\juerg{\,\ds}
\leq C \an{\norma{\ph }_{\L\infty V}^2 \norma{\dt \ph_1}_{\L2V}^2 }, $$
and as \eqref{reg:strong} holds, it turns out that 
\begin{align*}
	\norma{\an{g}}_{\L2 H} \leq	
	C \big ( 
	\norma{\ph_0}_{V}	
	+ \norma{w_{0}}_{V}	
		+ \norma{w'_{0}}_{}	
	+ \norma{\an{u}}_{\L2 H}	
	\big).
\end{align*}
Moreover, the initial value $y(0) =\a \an{v_0} + \b w_0 $ lies in $V$. 
Therefore, using parabolic regularity and the representation \an{given in} \eqref{pier6-9} (which holds 
as well), we easily infer that 
\begin{align}
	&\pico{\norma{w}_{\H1 H \cap \L\infty V\cap \L2 {\Hx2}}  + 
	\norma{\dt w}_{\H1 H \cap \L\infty V\cap \L2 {\Hx2}}}
	\non \\
	&\quad{}\leq C\big(\norma{\ph_0}_{V}+  \norma{w_0}_{V}+\norma{\an{v_0}}_V + \norma{\an{u}}_{\L2 H}\big). \label{pier12}
\end{align}}%

\noindent 
{\bf Third estimate:} \pier{First, we observe that \eqref{wf:cd:1} can be rewritten as 
\begin{align}
	\label{pier20}
	 \iO \dt\ph \, v = 
	- \iO \nabla \ph \cdot \nabla v
	+ \iO  	\pico{h}\hskip1pt v \quad\hbox{for every $v\in V$, a.e. in $(0,T)$}.
\end{align}
Here, recalling the notation in \eqref{not:diff:1}--\eqref{not:diff:3},
$\pico{h}$ is specified by
$$\pico{h} =  - \gamma (\ph_1)+ \gamma (\ph_2) - \frac2 {\theta_c} (\pi (\ph_1) - \pi (\ph_2))
+ \frac 1 {\theta_c^2} \big( \dt w  \, \pi (\ph_1) + \dt w_2 (\pi (\ph_1) - \pi (\ph_2))\big). $$
Now, in view of the regularity properties in \eqref{separation} and \eqref{reg:strong} \juerg{that
hold for both $(\ph_1, w_1)$ and  $(\ph_2, w_2)$,
we can check that every term of $\pico{h}$ belongs to $H^1(0,T;H)$ and that} 
\begin{align}
&\dt \pico{h}  = - (\gamma' (\ph_1) - \gamma' (\ph_2)) \dt \ph_1
         - \gamma' (\ph_2) \dt \ph 
- \frac2 {\theta_c} (\pi' (\ph_1) - \pi' (\ph_2))\dt \ph_1
- \frac2 {\theta_c} \pi' (\ph_2) \dt \ph 
\non \\ & \qquad {}
+ \frac 1 {\theta_c^2} \big( \partial_{tt} w  \, \pi (\ph_1) 
+ \dt w \, \pi'(\ph_1)\, \dt \ph_1 
+ \partial_{tt} w_2 (\pi (\ph_1) - \pi (\ph_2))\big)
\non \\ & \qquad {}
+ \frac 1 {\theta_c^2} \big( \dt w_2 (\pi' (\ph_1) - \pi' (\ph_2))\dt \ph_1
+ \dt w_2 \, \pi' (\ph_2) \dt \ph \big)\,.
\label{pier21} 
\end{align} 
Moreover, from \eqref{pier20} we can recover the expression of $\dt\ph (0) $, which is given by (cf.~\eqref{ass:exstrong:initialdata})
\begin{align} 
 \dt\ph (0) = \ph_{0,1}' - \ph_{0,2}' := &\ \Delta \ph_0 - 
(\gamma (\ph_{0,1})  - \gamma (\ph_{0,2}))
- \frac 2{\thc} ( \pi (\ph_{0,1})  - \pi (\ph_{0,2}))
\non \\ 
 &+ \frac 1 {\thc^2} \big( \an{v_0} \,\pi (\ph_{0,1}) 
 + \pico{v_{0,2}} (\pi (\ph_{0,1})  
  - \pi (\ph_{0,2}))\big)
 \label{pier22}  
\end{align}  
\juerg{and belongs to $H$, due to the assumptions on the} initial data. Therefore, since we also have that $\ph= \ph_1-\ph_2$ is in $\H1 V $,
\juerg{a comparison in \eqref{pier20} yields} that $\dt \ph \in \H1{V^*} $, and consequently we can differentiate \eqref{pier20} with respect to time and then test by 
$v= \dt \ph \in \L2V$. A subsequent integration leads to  
\begin{align} 
\frac 12 \norma{ \dt \ph (t)}^2 + \intQt |\nabla(\dt \ph)|^2 
= \frac 12 \norma{ \dt \ph (0)}^2 + \intQt \dt \pico{h} \, \dt\ph \,, 
\label{pier23}
\end{align} 
for all $t\in [0,T]$ (indeed, we also have $\dt \ph \in \C0H$\an{)}.   
Now, in view of \eqref{pier22} and \eqref{ass:exstrong:initialdata}, \eqref{strong:sep:initialdata}, it is straightforward to check that 
$$\frac 12 \norma{ \dt \ph (0)}^2 \leq 
C\big(\norma{\ph_0}_{W}^2 + \norma{\an{v_0}}_H^2 \big),$$
while, on account of the boundedness and Lipschitz continuity of $\gamma $ and $\pi$ in $[r_*, r^*]$
(cf.~\eqref{separation}), the \Holder\ and Young inequalities, and the continuous embedding 
$V\subset \Lx4$,  we can infer from \eqref{pier21} that 
\begin{align} 
  &\intQt \dt \pico{h} \, \dt\ph 
 \non \\ &\quad {}  
 \leq C \int_0^t \norma{\ph}_4 \norma{\dt \ph_1}_4 \norma{\dt \ph} \juerg{\,\ds} 
 + C \intQt |\dt\ph|^2  + C \intQt |\partial_{tt} w|^2 
 \non \\ &\qquad  {} 
  + C \int_0^t \norma{\dt w}_4 \norma{\dt \ph_1}_4 \norma{\dt \ph}  \juerg{\,\ds} 
  + C \int_0^t \norma{\partial_{tt} w_2 } \norma{\ph}_4 \norma{\dt \ph}_4 \juerg{\,\ds} 
 \non \\ &\quad {} 
  \leq C \big( \norma{\ph}_{\L\infty V} + \norma{\dt w }_{\L\infty V} \big)
  \norma{\dt \ph_1}_{\L2V} 
  \norma{\dt \ph }_{\L2H}  + C\norma{\dt \ph }_{\L2H}^2
  \non \\ &\qquad {}  
    + C \norma{\partial_{tt} w}_{\L2H}^2 
  \big(1 +\norma{\ph}_{\L\infty V}^2 \big)
  + \frac 12 \intQt \big(|\dt \ph (t)|^2 +  |\nabla(\dt \ph)|^2 \big). \non
  \end{align} 
Then, by virtue of \eqref{pier10} and \eqref{pier12}, combining the last two inequalities with \eqref{pier23} plainly leads to the estimate} 
 \begin{align}\label{pier26}
	\pier{\norma{\ph}_{\W{1,\infty} H \cap \H1 {V}}
	\leq C \big ( 
	\norma{\ph_0}_{W}	
	+ \norma{w_{0}}_{V}	
		+ \norma{\pico{v}_{0}}_{V}	
	+ \norma{\an{u}}_{\L2 H}	
	\big).}
\end{align}

\noindent
{\bf \pier{Fourth} estimate:}
\pier{Now, from \eqref{pier20}, that reproduces~\eqref{wf:cd:1}, and the regularity of solutions we deduce that 
\begin{align*}
	-\Delta \ph =  \pico{h} - \dt \ph 
	\quad\hbox{a.e. in }\, Q,	
\end{align*}
with the \rhs\ that is under control in $\L\infty H$. Then, by elliptic regularity we easily derive the estimate 
\begin{align}
	\norma{\ph}_{\L\infty {\Hx2}} 
		\leq 	
	C \big ( 
	\norma{\ph_0}_{W}	
	+ \norma{w_{0}}_{V}	
		+ \norma{\pico{v}_{0}}_{V}	
	+ \norma{\an{u}}_{\L2 H}	
	\big).\label{pier11}
\end{align}}%

Therefore, upon collecting \eqref{pier10}, \eqref{pier12}, \eqref{pier26}, \an{and} \eqref{pier11}, we obtain \eqref{cd:strong} and conclude the proof of Theorem~\ref{THM:UQ:STRONG}.
\end{proof}

\section{Optimal control theory}
\label{SEC:OC}
\setcounter{equation}{0}

In this section, we aim at solving an optimal control problem whose governing state equation is given
 by the system \Sys\ analyzed in the previous section.
We seek optimal controls in the form of a distributed heat source, represented by $\an{u}$ in \eqref{sys2}, and an initial temperature, 
which corresponds to $\an{v_0}$ in \eqref{sys4}. 
\an{A}s we aim at covering the cases of polynomial and regular singular potentials, including, in particular, \eqref{Freg} and \eqref{Flog}, we are from now on restricting ourselves to the framework of strong solutions (cf. Theorems \ref{THM:EX:STRONG} and \ref{THM:UQ:STRONG}). 

The control problem \juerg{under investigation reads as follows}:

\vspace{3mm}\noindent
\CP \quad Minimize the cost functional
\begin{align}\label{cost}
		\notag
		 \J(u, v_0, \vp, w)  &: = \frac {k_1} 2 \norma{\vp - \vp_Q}_{L^2(Q)}^2
		+ \frac {k_2} 2 \norma{\vp(T) - \vp_\Omega}^2
		+ \frac {k_3} 2  \norma{ w - w_Q}_{L^2(Q)}^2 
		\\ & \quad 	\notag
		+ \frac {k_4} 2  \norma{ w(T) - w_\Omega}^2
		+ \frac {k_5} 2 \norma{\dt w - w_Q'}_{L^2(Q)}^2 
		+ \frac {k_6} 2  \norma{\dt w(T) - w'_\Omega}^2
		\\ & \quad 
		+ \frac {\nu_1} 2 \norma{u}_{L^2(Q)}^2 + \frac {\nu_2} 2 \norma{v_0}_V^2 
\end{align}
subject to the state system \an{\Sys\ and}
to the control constraint
\begin{align*}
	(u,v_0) \in \Uad,
\end{align*}
where $\UU := L^\infty(Q) \times (V\cap \Lx{\infty} )$ and \an{the set of {\it admissible controls} is}
\begin{align}
	&\Uad :=
		\big\{ 			
			(u,v_0) \in \UU :  u_* \leq u \leq u^* \,\,\aeQ, \non \\
			&\qquad\qquad\qquad\qquad v_* \leq v_0 \leq v^* \,\,\aeO, \quad \norma{v_0}_V \leq M
		\big\}.
		\label{Uad}
\end{align}
Above, the symbols $k_1, ..., k_6$ and $\nu_1,\nu_2 $ denote some nonnegative constants which are not all zero, while \juerg{\,$\vp_Q, w_Q, w_Q'\in L^2(Q)$ \,and $\,\vp_\Omega,w_\Omega,w_\Omega'\in L^2(\Omega)$  denote some prescribed targets}.
As for the set of admissible controls $\Uad$, we 
assume that $u_* $ and $ u^*$ are prescribed functions in $ L^\infty(Q)$; moreover,  $ v_*$ and $v^*$ 
are given in $L^\infty(\Omega)$, and  $\,M>0\,$  is a fixed constant such that  
\begin{align}
\non 
	\Uad \, \hbox{ is a nonempty, closed and convex subset of the control space } \,  \UU.
\end{align}
Note that closedness and convexity can be easily verified from \eqref{Uad}. 
Furthermore, we can select a value $R>0$ big enough such that the open ball
\begin{align}\label{UR}
	\UR :=
		\big\{ 			
			(u,v_0) \in \UU:  \norma{(u,v_0)}_\UU < R
		\big\} \, \hbox{ contains } \, \Uad.
\end{align}

Let us remark that from a physical viewpoint it is more relevant investigating the evolution of $\dt w$ 
instead that of $w$, as the first one denotes the temperature of the system.
This is the reason why the terms in \eqref{cost} related to $k_5$ and $k_6$ are more significant than the ones associated with $k_3$ and $k_4$; nonetheless, we believe that those less physical terms are still worth 
\juerg{considering} from a mathematical viewpoint through the way \juerg{in which} they appear in the 
adjoint system (cf. system \eqref{adsys1}--\eqref{adsys4}). \pier{Also, note that the quantities 
$ v_*$ and $v^*$ appearing in \eqref{Uad} represent threshold values for the initial 
\juerg{temperature distribution}\, $v_0$, while the condition
$\norma{v_0}_V \leq M$ prevents extremely large variations for this distribution.}  

By virtue of \pier{Theorems~\ref{THM:EX:WEAK}--\ref{THM:UQ:STRONG}}, the
{\it control-to-state operator}
\begin{align*}
	\S: \UU_R \subset \UU  \to \Y  , \quad \S:(u,v_0) \mapsto (\ph,w),
\end{align*}
is well-defined as a mapping from $\UU$ into the solution space $\Y$, with the latter being defined by (cf. Theorem \ref{THM:EX:STRONG})
\begin{align}
	\non
	\Y & := \bigl({\W{1,\infty} H \cap \H1 V} \cap \L\infty{W}\bigr)
	\\ & \qquad \times%
	\bigl(\H2 H \cap \W{1,\infty} V \cap \H1 {W}\bigr).
	\label{sol:space}
\end{align}
Moreover, \pico{we also set}
\begin{align*}
	\non
	\X & := \bigl(\H1 H \cap \L\infty V \cap \L2 W\bigr)
	\\ & \qquad \times
	\bigl(\H2 H \cap \W{1,\infty} V \cap \H1 {W}\bigr)
\end{align*}
\pico{and observe that $\Y\subset \X$ with continuous embedding.}
Then, the solution operator allows us to define the {\it reduced cost functional} as follows:
\begin{align}\label{red:cost}
	\Jred: \UU \to \erre, \quad \Jred(u,v_0) := \J \big(u,v_0, \S(u,v_0)\big).
\end{align}
Moreover, notice that Theorems \ref{THM:UQ:WEAK} and \ref{THM:UQ:STRONG} already ensure that the solution operator $\S$
is \Lip\ continuous \pier{in $\UU_R$ when viewed as a mapping from $L^2(Q)\times V$} into the space $\an{\Y}$. Namely, 
for arbitrary controls $(u_i,v_{\pier{{0,i}}}) \in \UU_R$, \pier{$i=1,2$}, the stability estimate \eqref{cd:strong} yields that
\begin{align*}
	\norma{\S(u_1,v_{\pier{{0,1}}})-\S(u_2,v_{\pier{{0,2}}})}_\an{\Y}
	\leq C 	\big( \norma{u_1-u_2}_{\L2 H}
	+\norma{v_{\pier{{0,1}}}- v_{\pier{{0,2}}}}_{V}\big).
\end{align*}

For the control problem, some additional assumptions are in order:
\begin{enumerate}[label={\bf C\arabic{*}}, ref={\bf C\arabic{*}}]
\item \label{ass:reg:non:control}
 $\hat \gamma\in \pier{C^3(r_-, r_+)}$ and $\hat \pi \in C^3(\erre)$.
 \item \label{ass:control:const}
	$k_1,k_2,k_3,k_4,k_5,k_6,\nu_1,\nu_2$ are nonnegative constants, not all zero. 
\item \label{ass:control:target}
	The target functions fulfill $\vp_Q ,w_Q,  w_Q' \in L^2(Q)$, $\vp_\Omega, {w_\Omega}\in H$, and $w_\Omega' \in V$.
\item \label{ass:control:Uad}
	The functions $u_* , u^*$ belong to $ L^\infty(Q)$ with $u_* \leq u^* \,\, 
	\aeQ$, and \pier{$ v_*,v^*$ are fixed in $L^\infty(\Omega)$ such that $ v_*\leq v^* \,\,\aeO$. 
	\juerg{Moreover, $M>0$, and} the set $\Uad$ defined by \eqref{Uad} is nonempty.}
\end{enumerate}

The first result we address concerns the existence of an optimal strategy, that is of an optimal control pair.
\begin{theorem}\label{THM:CON:EX}
Suppose that \ref{ass:1:constants}--\ref{ass:3:pi}, \ref{ass:strong:gamma:1}--\ref{ass:4:pidiff}, \an{\ref{ass:control:const}}--\ref{ass:control:Uad} hold
\pier{in addition to the assumptions~\eqref{ass:exweak:initialdata}, \eqref{ass:exstrong:initialdata}, \eqref{strong:sep:initialdata}
on $\ph_0, w_0$.} Then the minimization problem \CP\ admits \pier{a solution},
that is, there exists at least one optimal pair $\opt \in \Uad$ such that
\begin{align*}
		\Jred (\bu,\bv) \leq \Jred (u,v_0) 
		\quad \forall (u,v_0) \in \Uad. 
\end{align*}
\end{theorem}
\begin{proof}[Proof of Theorem \ref{THM:CON:EX}]
\pier{The existence of a minimizer $\opt$ plainly follows from applying the direct method of the 
calculus of variations. \pico{In fact, we can pick} a minimizing sequence $\an{\{(u_n,v_{0,n})\}_n}\pico{{}\subset{}} \Uad $ for the 
functional $\Jred $, and let\an{, for every $n \in \mathbb{N}$,} $(\ph_n, w_n)= \S (u_n,v_{0,n}) $ denote the corresponding 
strong solution to the system \pico{\eqref{sys1}--\eqref{sys4}}. Then, due to \eqref{Uad}, by 
compactness it turns out that there exist a subsequence, still denoted by $\an{\{(u_n,v_{0,n})\}_n}$,
and a pair $(\bu,\bv) \in \Uad$ such that 
\begin{align} 
 u_n &\to \bu \quad \hbox{weakly star in } L^\infty (Q), \non \\
 v_{0,n} &\to \bv \quad \hbox{weakly star in } V \cap L^\infty (\Omega), \non 
\end{align} 
as $n\nearrow \infty$. Correspondingly, in view of Theorem~\ref{THM:EX:STRONG},
 and taking advantage of \cite[Sect.~8, Cor.~4]{Simon}, it turns out that there is a 
pair $(\bph,\bw) $ satisfying
\begin{align} 
\ph_n \to \bph \quad &\hbox{weakly star in }  \W{1,\infty} H \cap \H1 V \cap \L\infty{W}
	\non \\ &\hbox{and strongly in } C^0(\overline Q) 
, \label{pier13} \\
 w_n \to \bw \quad &\hbox{weakly star in }  \H2 H \cap \W{1,\infty} V \cap \H1 {W}
	\non \\ & \hbox{and strongly in } \C1H \cap \H1V, 
	 \label{pier14} 
\end{align}
in principle for another subsequence. Indeed, as for \eqref{pier13} note that $W \subset C^0(\overline \Omega)$ with compact embedding.
At this point, it is a standard matter to check that 
\juerg{passage} to the limit as $n\nearrow \infty $ in the system \an{\Sys},
written for $\{\ph_n, w_n, u_n, v_{0,n}\}$, leads to the same system written for \an{the limits}
$\{\bph, \bw, \bu, \bv\}$. Then, taking into account Theorem~\ref{THM:UQ:STRONG}
as well, we infer that  $(\bph, \bw)= \S (\bu,\bv) $, \eqref{pier13} and \eqref{pier14}
hold for the selected subsequence, and, by the lower semicontinuity of norms, 
\begin{align*}
		\Jred (\bu,\bv) \leq \liminf_{n\to \infty}  \, \Jred (u_n,v_{0,n} ).
\end{align*}
\juerg{Hence,} $\opt $ is a global minimizer for $\Jred$, as $ \Jred (u_n,v_{0,n} )$ converges 
to the infimum of $\Jred$. The assertion is thus proved.}
\end{proof}

We are now interested in finding optimality conditions that every minimizer \an{has to} satisfy.
\pier{To this end,} recall the reduced form \eqref{red:cost} and the fact that $\Uad$ is a nonempty, closed, and convex subset of the control space $\UU$. Standard results of convex analysis (see, e.g., \cite{Trol}) \pier{entail} the first-order necessary condition for $\Jred$ at every minimizer $\opt$ in terms of a suitable variational inequality of the form
\begin{align}
	\label{formal:foc}
	D\Jred \opt (u- \bu, v_0 - \bv) 
	\geq 0 \quad  \forall(u,v_0) \in \Uad,
\end{align}
where $D\Jred $ stands for the derivative of the reduced cost functional
in a proper mathematical sense (cf. Theorem~\ref{THM:CON:FRE}).
The quadratic structure of $\J$ directly yields its \Frechet\ differentiability, so that it suffices
to show the differentiability of the solution operator $\S$ in order to derive the first-order necessary conditions from \eqref{formal:foc} by means of the chain rule.

\juerg{For this purpose, we} fix a control pair $(\bu,\bv) \in \UU_R$ 
with corresponding state $ (\bvp, \bw) = \S(\bu,\bv)$. We introduce the linearized system \pier{to \an{\Sys}}, which reads, for every $(h, h^0) \in L^2(Q) \times V$, 
\begin{alignat}{2}
\label{linsys1}
&\dt\xi-\Delta\xi + \gamma'(\bvp)\xi + \tfrac 2 {\thc} \pi'(\bvp)\xi 
	- \tfrac 1{\thc^2}{\dt \et}\,\pi(\bvp)
	- \tfrac 1{\thc^2}{\bdtw}\,\pi'(\bvp)\xi = 0 \quad &&\mbox{a.e. in }Q,\\
\label{linsys2}
&\ddt \et -\alpha \Delta (\dt \et)-\beta \Delta \et + \pi' (\bvp)\xi \dt \bvp+ \pi (\bvp) \dt \xi =h \quad && \mbox{a.e. in }Q,\\
\label{linsys3}
&\dn \xi=\dn (\alpha \dt \et + \beta \et)=0 \quad &&\mbox{a.e. on }\Sigma,\\
\label{linsys4}
&\xi(0)=0,\quad \et(0)=0,\quad \dt \et(0)=h^0\quad &&\mbox{a.e. in }\Omega.
\end{alignat}
\Accorpa\Lin {linsys1} {linsys4}
Its well-posedness is stated in the following result.

\begin{theorem}\label{THM:WP:LIN}
	Assume that \ref{ass:1:constants}--\ref{ass:3:pi} 
	and \ref{ass:strong:gamma:1}--\ref{ass:4:pidiff} are fulfilled \pier{in addition to the assumptions~\eqref{ass:exweak:initialdata}, \eqref{ass:exstrong:initialdata}, \eqref{strong:sep:initialdata}
on $\ph_0, w_0$.} Let \juerg{$(\bu,\bv)\in\UU_R$ be given and} $ (\bvp, \bw) = \S(\bu,\bv)$.	
	Then \juerg{the linearized system \Lin\ 	has for every $(h, h^0) \in L^2(Q) \times V$ a} 
	unique solution $(\xi,\et)\in \X$.
\end{theorem}

\begin{proof}[Proof of Theorem \ref{THM:WP:LIN}]
Since the problem is linear, we can prove existence and, at the same time, uniqueness, by performing suitable estimates on the solution $(\xi,\et)\in \X$  in terms of the data  $(h, h^0) \in L^2(Q) \times V$, with linear dependence. \juerg{As in the case of the state problem, we here} avoid to implement a Faedo--Galerkin scheme 
\juerg{and argue}  directly on the linearized problem.

\noindent
{\bf First estimate:}
\juerg{We first add $\xi $ to both sides of \eqref{linsys1} and then} test \eqref{linsys1} by $\thc^2 \dt \xi$ and \eqref{linsys2} by $\dt \et$. Next, we sum up the resulting equalities and integrate by parts to infer that a cancellation occurs\juerg{, obtaining the identity}
\begin{align*}
	& \thc^2 \I2 {\dt \xi}
	+ \frac {\thc^2} 2 \IO2 {\xi}_V
	+ \frac 12\IO2 {\dt \et}
	+ \a \I2 {\nabla (\dt\et)}
	+ \frac \b2 \IO2 {\nabla \eta}
	\\ & \quad 
	 =  \frac 12 \norma{h^0}^2
	 + \intQt \big(\thc^2 - \thc^2\gamma'(\bvp) -  2 {\thc}\pi'(\bvp) ) \xi \, \dt \xi
	\\ & \qquad 
	+ \intQt {\bdtw}\,\pi'(\bvp)\xi \, \dt \xi
	- \intQt \pi' (\bvp) \dt \bvp\, \xi \, \dt \et 
	+ \intQt h \,\dt\et= : \sum_{i=1}^{5} \II_i.
\end{align*}
Since $(\bvp ,\bw) $ is a strong solution to \Sys,
we deduce from \eqref{separation}--\eqref{reg:strong} that 
\juerg{$\gamma'(\bvp), \linebreak \pi'(\bvp), \, \dt\bw \in L^\infty(Q)$ and 
$\dt\bph \in {\L\infty H \cap \L2 V}$. We thus have that}
\begin{align*}
	\II_2 + \II_3
	\leq 
	\frac {\thc^2} 2  \I2 {\dt \xi}
	+ C \intQt  | \xi|^2 \,,
\end{align*}
and, with the help of the continuous embedding $V\subset \Lx4$, 
\begin{align*}
	\II_4
	\leq 
	C \int_0^t \norma{ \dt\bph }_4 	\norma{\xi}_4 \an{\norma{ \dt \et}} \juerg{\,\ds} 
		\leq 
	C \int_0^t \norma{ \dt\bph }_V^2 	\norma{\xi}_V^2 \juerg{\,\ds} 
	+ C \intQt |\dt \et|^2 \,,
\end{align*}
where the function $t \mapsto  \norma{ \dt\bph (t) }_V^2$ belongs to $L^1(0,T)$.
As for the last term, we simply employ Young's inequality and obtain
\begin{align*}
	\II_5 \leq 
	C \intQt( |{h}|^2 + |\dt \et|^2).
\end{align*}
We collect the above estimates and apply Gronwall's lemma. Then, observing that 
(cf.~\eqref{linsys4}) \,\,$\norma{\eta(t)}_V^2 \leq T \intQt \norma{\dt \eta}_V^2 $ for $t\in [0,T]$, by the 
H\"older inequality, we can conclude that 
\begin{align}
	\norma{\xi}_{\H1 H \cap \L\infty V }
	+\norma{\et}_{\W{1,\infty} H \cap \H1  V}
	\leq C (\norma{h}_{L^2(Q)} + \norma{h^0}).
	\label{pier15} 
\end{align}

\noindent
{\bf Second estimate:}
Next, \pier{\eqref{pier15} (in particular, the boundedness of
$\Vert\dt \xi \Vert_{\L2 H }$) and a comparison of terms in \eqref{linsys1} easily produce that} 
\begin{align*}
	\norma{\Delta\xi}_{\L2 H}
	\leq C (\norma{h}_{L^2(Q)} + \norma{h^0}),
\end{align*}
so that elliptic regularity entails that 
\begin{align}
	\norma{\xi}_{\L2 W}
	\leq C (\norma{h}_{L^2(Q)} + \norma{h^0}).
	\label{pier16}
\end{align}

\noindent
{\bf Third estimate:}
As done in the third estimate of \pier{Theorem~\ref{THM:EX:STRONG}}, we add to both sides of \eqref{linsys2} the term $\frac \b\a \dt\et$ and rewrite it as a parabolic \pier{equation} in terms of the new variable $y:=\a \dt \et + \b \et$. \pier{Precisely, we deduce that
\begin{align*}
  \begin{cases}
 	\frac 1\a\dt y -\Delta y = g \quad & \text{in}\,\, Q,
 	\\ 
 	\dn y = 0 \quad & \text{on}\,\, \Sigma,
 	\\ 
 	\pier{y(0) = \alpha h^0} \quad & \text{in}\,\, \Omega,
  \end{cases}
\end{align*}
with \,$g :=  -\pi' (\bvp)\xi \dt \bvp -  \pi (\bvp) \dt \xi  + h + \frac \b\a \dt \et$, here.
Due to \eqref{pier15}, we have that the norm of $g$  in $\L2H$ is under control. Besides, we are assuming that $h^0 \in V$, so that parabolic regularity theory entails that}
\begin{align}
	\norma{\a \dt \et+ \b \et}_{\H1 H \cap \L\infty V \cap \L2 W} 
	\leq C (\norma{h}_{L^2(Q)} + \norma{h^0}_V).
	\label{pier17}
\end{align}
Now, \pier{arguing as in \eqref{pier6-9}, 
it follows that \eqref{pier17} implies the same estimate for $\et $ and $\dt \et$, whence}
\begin{align}
	\norma{\et}_{\H2 H \cap \W{1,\infty} V \cap \H1 W} 
	\leq C (\norma{h}_{L^2(Q)} + \norma{h^0}_V).
	\label{pier18}
\end{align}
Then, by collecting \eqref{pier15}, \eqref{pier16}, and \eqref{pier18}, we end the proof.
\end{proof}

After \pier{proving} Theorem~\ref{THM:WP:LIN}, we are in a position to show that the control-to-state operator $\S$ is \Frechet\ differentiable as a mapping between suitable Banach spaces.
Here is the \pier{related} result.
\begin{theorem}\label{THM:CON:FRE}
Suppose that the conditions \ref{ass:1:constants}--\ref{ass:3:pi}, \ref{ass:strong:gamma:1}--\ref{ass:4:pidiff}, and~\ref{ass:reg:non:control} are 
fulfilled.
Moreover, let the initial data \pier{$\ph_0$ and $w_0$ satisfy \eqref{ass:exweak:initialdata}, \eqref{ass:exstrong:initialdata}, \eqref{strong:sep:initialdata}, and let $(\bu, \bv) \in \UU_R$  with $(\bvp,\bw)=\S(\bu, \bv)$. 
Then the solution operator $\S$ is \Frechet\ differentiable at $(\bu, \bv)$ as a mapping from $\UU$ into 
$\X$. Moreover, for every $\bh : =(h,h^0)\in \UU$, the \Frechet\ derivative $D \S(\bu,\bv)\in 
{\cal L}(\cal U,\X)$ is given by the identity $D\S(\bu, \bv)({\bh})=(\xi,\et)$, where 
$(\xi,\et)$ is the unique solution to the linearized system \eqref{linsys1}--\eqref{linsys1} associated with $\bh$.}
\end{theorem}

\begin{proof}[Proof of Theorem \ref{THM:CON:FRE}]
Since $\UR$ is open, provided that we consider small $\eps$-perturbations in the $\UU$-norm, we surely have that $(\bu+ h, \bv + h^0)\in\UR$ as well, that is, there exists \juerg{some} $\eps >0$ such that
\begin{align*}
	(\bu+ h, \bv + h^0)\in\UR \quad \forall \, \bh \in \UU \quad \text{such that} \quad \norma{\bh}_\UU \leq \eps.
\end{align*}
For the rest of the proof, we agree that this condition is met by all of the appearing increments $\bh$.

\juerg{We} claim that $D \S \opt (\bh) = (\xi,\et)$, with $(\xi,\et)$ being the unique solution to the linearized system \Lin. \juerg{We prove this claim directly by showing that} 
\begin{align}
	\S (\bu + h,\bv +h^0) = \S (\bu,\bv) + \pier{(\xi,\et)} + o(\norma{\bh}_{\UU})
   \quad \hbox{in $\cal X$} \quad      \hbox{as} \ \  \norma{\bh}_{\UU}\to 0.
   \label{Fre:formal}
\end{align}
Upon setting 
\begin{align}\label{not:fre}
	 (\bvph, \bwh)= \S(\bu+h,\bv+h^0 ), 
	 \quad 
	 \ps:= \bvph - \bvp - \xi, 
	\quad 
	\pier{z} := \bwh - \bw - \et,
\end{align}
the condition \eqref{Fre:formal} becomes
\begin{align}\label{fre:final}
		\norma{(\psi, \pier{z})}_\X = o(\norma{\bh}_{\UU})
   \quad \hbox{as} \ \  \norma{\bh}_{\UU}\to 0,
\end{align}
which is the identity we are going to prove. Accounting for the notation in \eqref{not:fre}, we infer that the variables $\psi$ and $\pier{z}$ solve the initial-boundary value problem
\begin{align}
\label{fresys1}
	&\dt\psi-\Delta\psi  + \Lambda_1 \pier{{}- \an{\tfrac 1{\thc^2}} \pi(\bvp)  \dt \pier{z}}=0  \quad&&\mbox{in }\,Q\,,\\
\label{fresys2}
	&\ddt \pier{z} -\alpha \Delta (\dt \pier{z})-\beta \Delta \pier{z} + \Lambda_2 {+\pi (\bvp) \dt \psi} =0 \quad&&\mbox{in }\,Q\,,\\
\label{fresys3}
	&\dn \psi=\dn (\alpha \dt \pier{z} + \beta \pier{z})=0 \quad&&\mbox{on }\,\Sigma,\\
\label{fresys4}
	&\psi(0)=0,\quad \pier{z}(0)=0,\quad \dt \pier{z}(0)=0\,\quad &&\mbox{in }\,\Omega,
\end{align}
\Accorpa\Fre {fresys1} {fresys4}
where \pier{the terms $\Lambda_1$ and $\Lambda_2$} are defined by
\begin{align*}
	& \Lambda_1  = [\gamma(\bvph) - \gamma(\bvp) -\gamma'(\bvp)\xi]
	+ \an{\tfrac 2{\thc}}  [\pi(\bvph) - \pi(\bvp) -\pi'(\bvp)\xi]
	\\ & \qquad  
	\an{-\tfrac 1{\thc^2}}
	\Big(
	(\pi(\bvph)-\pi(\bvp))(\dt \bwh-\dt \bw)
	+ \dt\bw \,[\pi(\bvph) - \pi(\bvp) -\pi'(\bvp)\xi]\Big)
	,
	\\[2mm]
	& \Lambda_2 = 
	(\pi(\bvph)-\pi(\bvp))(\dt \bvph-\dt \bvp)
	+ \dt\bvp \,[\pi(\bvph) - \pi(\bvp) -\pi'(\bvp)\xi].
\end{align*}

Before moving on, let us recall that the continuous dependence result in Theorem \ref{THM:UQ:STRONG}, applied to the solutions $ (\bvph, \bwh)$
and  $(\bvp, \bw)$, yields that 
\begin{align}\label{cd:fre}
	\non
	& 
	\an{ \norma{\bvph- \bvp}_{\W{1,\infty} H \cap \H1 V \cap \L\infty {W}}}
	+ 
	\norma{\bwh -\bw}_{ \H2 H \cap \W{1,\infty} V \cap \H1 {\Hx2}}
	\\ & \quad 
	\leq \pier{K_4 \big( \norma{h}_{L^2(Q)} + \norma{h^0}_V \big)}.
\end{align}
Besides,  $ (\bvph, \bwh)$ and  $(\bvp, \bw)$, as strong solution to \pier{\an{\Sys},}  
satisfy \pier{\eqref{separation} and \eqref{reg:strong}}. 
Moreover, we recall Taylor's formula with integral remainder:
let $g :\erre \to \erre$ be a differentiable function with \Lip\ continuous derivative $g'$. Then, \pier{for $\overline x \in \erre$} it holds that
\begin{align}
	\label{taylor:formal}
	g(x)= g(\overline x) + g'(\overline x) (x - \overline x) 
	+  (x - \overline x)^2  \int_{0}^1 g''( \overline x + s (x - \overline x))(1-s) \, \ds\an{,	\quad x \in \erre.}
\end{align}
An application of \eqref{taylor:formal} to $\pi$ and $\gamma$ yields that
\begin{align}
	\label{tay:gamma}
	\gamma(\bvph) - \gamma(\bvp) - \gamma'(\bvp)\xi &= \gamma'(\bvp) \psi + R_{\gamma}^\bh \,(\bvph-\bvp)^2,
	\\ 
	\label{tay:pi}
	\pi(\bvph)-\pi(\bvp)-\pi'(\bvp)\xi &= \pi'(\bvp) \pier{\psi}  + R_{\pi}^\bh \,(\bvph-\bvp)^2,
\end{align}
with the remainders
\begin{align*}
	R_{\gamma}^\bh:= \int_0^1 \gamma''(\bvp+s (\bvph-\bvp)) (1-s){\,\ds}, 
	\quad
	R_{\pi}^\bh:= \int_0^1 \pi''(\bvp+s (\bvph-\bvp)) (1-s){\,\ds}.
\end{align*}
Due to assumptions \ref{ass:reg:non:control}, it directly follows that 
\begin{align}\label{unif:rem}
	\norma{R_{\gamma}^\bh}_{L^\infty(Q)}
	+\norma{R_{\pi}^\bh}_{L^\infty(Q)}
	\leq C.
\end{align}

We now \pier{prove} some estimates that will imply \eqref{fre:final}.

\noindent
{\bf First estimate:}
\pier{Add $\psi$ to both sides of \eqref{fresys1} and test it by $\thc^2 \dt \psi$; then, test \eqref{fresys2} by $\dt \pier{z}$ and sum up} the resulting equalities.
After integration by parts, we obtain that a cancellation occurs and that
\begin{align*}
	& \pier{\thc^2\I2 {\dt \psi}
	+ \frac{\thc^2}2 \pier{\IO2 {\psi}_V}}
	+ \frac 12\IO2 {\dt \pier{z}}
	+ \a \I2 {\nabla (\dt\pier{z})}
	+ \frac \b2 \IO2 {\nabla \pier{z}}
	\\ & \quad 
	 = {}
	 \pier{{}\thc^2 \intQt (\psi - \Lambda_1) \dt \psi{}}
	 - \intQt \Lambda_2 \, \dt \pier{z} .
\end{align*}
\pier{The \pier{first} term on the \rhs\ can be controlled by employing 
Taylor's formulae \eqref{tay:gamma}--\eqref{tay:pi}, the uniform bounds 
\eqref{separation}--\eqref{reg:strong} for $(\bvp, \bw)$, the
Young and \Holder\ inequalities, the stability estimates \eqref{cd:fre}, 
\eqref{unif:rem}, \juerg{and} the continuous embedding~$V \subset \Lx4$. We} infer that 
\begin{align}
	&  \pier{\thc^2 \intQt (\psi - \Lambda_1) \dt \psi }	
	\non\\ 
	& \quad 	\leq 
	\thc^2 \intQt \big| \pier{{}\psi -  \gamma'(\bvp) \psi - 
	R_{\gamma}^\bh (\bvph-\bvp)^2{}}\big|  | \dt\psi| 
	+ 2 \thc \intQt \big|\pi'(\bvp) \psi  + R_{\pi}^\bh (\bvph-\bvp)^2\big|  | \dt\psi|
	\non\\
	& \qquad 
	+ \intQt |\pi(\bvph)-\pi(\bvp) | |\dt \bwh-\dt \bw| |\dt\psi|
	+ \intQt |\dt\bw|  \big|\pi'(\bvp) \psi  + R_{\pi}^\bh (\bvph-\bvp)^2\big|  | \dt\psi|
	\non\\ 
	& \quad 
	\leq 
	\d \I2 {\dt \psi}
	+ \cd \intQt |\psi|^2
	+ \cd \iot \norma{\bvph-\bvp}^4_4 \juerg{\,\ds}
	\non\\                        
	& \qquad 
	+ \cd \iot \norma{\bvph-\bvp}^2_4\norma{\dt \bwh-\dt \bw}^2_4 \juerg{\,\ds}
    \non\\ & \quad 
	\leq 
	\d \I2 {\dt \psi}
	+ \cd \intQt |\psi|^2
	+ \cd \iot \norma{\bvph-\bvp}_V^2 
	 \big( \norma{\bvph-\bvp}_V^2+ \norma{\dt \bwh-\dt \bw}^2_V \big)\juerg{\,\ds}
		\non\\ & \quad 
	\leq 
	\d \I2 {\dt \psi}
	+ \cd \intQt |\psi|^2
	+ \cd \pier{\big( \norma{h}_{L^2(Q)}^4 + \norma{h^0}_V^4 \big)},
	\label{lambda1} 
\end{align}
for a positive $\d$ yet to be chosen. Similar arguments allow us to bound 
the \pier{second} term on the \rhs, concluding that 
\begin{align}
	&{} - \intQt \Lambda_2 \, \dt \pier{z} \non \\
	&\quad \pier{{}\leq \intQt
	|\pi(\bvph)-\pi(\bvp)|  |\dt \bvph-\dt \bvp | |\dt \pier{z}|
	+ \intQt \big|\dt\bvp | \big| \pi'(\bvp) \pier{\psi}  + R_{\pi}^\bh 
	\,(\bvph-\bvp)^2 \big|
	|\dt \pier{z}| }\non \\
	&\quad \pier{{}\leq 
	C \iot \norma{\bvph-\bvp}_4\norma{\dt \bvph-\dt \bvp}_4 \norma{\dt \pier{z}}\juerg{\,\ds}
	+ C \iot  \norma{\dt \bvp}_6 \norma{\dt \pier{z}} \big(
	\norma{\psi}_3 + \norma{\bvph-\bvp}^2_6	\big)}\juerg{\,\ds}
	\non \\
	&\quad	\pier{{} \leq 
    C \iot \big(1 + \norma{\dt \bvp}_V^2\big)  \norma{\dt \pier{z}}^2	\juerg{\,\ds} 	
    + C \norma{\bvph-\bvp}_{\L\infty V}^2 \iot \norma{\dt \bvph-\dt \bvp}_V^2}\juerg{\,\ds}
    \non \\ &\qquad	
	\pier{{}
	+ C \iot \norma{\psi}_V^2 \juerg{\,\ds} + C \iot \norma{\bvph-\bvp}^4_V }\juerg{\,\ds}
	\non \\ & \quad \pier{{}\leq 
	C \iot \big(1+ \norma{\dt \bvp}^2_V \big) \norma{\dt \pier{z}}^2	\juerg{\,\ds}
	+ C \iot \norma{\psi}_V^2}\juerg{\,\ds}
	\pier{{} + C \pier{\big( \norma{h}_{L^2(Q)}^4 + \norma{h^0}_V^4 \big)},}
	\label{lambda2} 
\end{align}
\pier{where we notice that the function $t \mapsto \an{(1+ \norma{\dt \bvp}^2_V )}$ is in $ L^1(0,T)$,
 due to \eqref{reg:strong}.
Upon choosing $0<\d< \theta_c^2$,} Gronwall's lemma yields that
\begin{align}
	\norma{\psi}_{\H1 H \cap \L\infty V }
	+\norma{\pier{z}}_{\W{1,\infty} H \cap \H1  V}
	\leq C \pier{\big( \norma{h}_{L^2(Q)}^2 + \norma{h^0}_V^2 \big)}.
	\label{pier19}
\end{align}

\noindent
{\bf Second estimate:}
\pier{A closer inspection of the estimate in \eqref{lambda1}, along with the bound \eqref{pier19},
 shows that $\norma{\Lambda_1}_{\L2 H}$ is bounded as well by an analogous term. Then, a comparison 
argument in \eqref{fresys1} directly leads to}
\begin{align*}
	\norma{\Delta \psi}_{\L2 H} 
	\leq C \pier{\big( \norma{h}_{L^2(Q)}^2 + \norma{h^0}_V^2 \big)},
\end{align*}
so that \pier{\eqref{pier19} and elliptic regularity yield} \juerg{that}
\begin{align}
\label{pier24}
	\norma{\psi}_{\L2 W} 
	\leq C \pier{\big( \norma{h}_{L^2(Q)}^2 + \norma{h^0}_V^2 \big).}
\end{align}

\noindent
{\bf Third estimate:}
\pier{Repeating the argument employed in the third estimate of the proof of 
Theorem~\ref{THM:EX:STRONG} (cf., in particular, \eqref{parab:sys}), we can
in view of \eqref{fresys2}--\eqref{fresys4} state a parabolic system in the variable $y= \a \dt z + \b z$, with source term $ \tfrac \b\a \dt z - \Lambda_2  - \pi (\bvp) \dt \psi $ and null initial value. With the help of \eqref{pier19}, it is not difficult to verify that 
\begin{align*}
	\big\| \tfrac \b\a \dt z - \Lambda_2  - \pi (\bvp) \dt \psi \big\|_{\L2 H} 
	\leq C \pier{\big( \norma{h}_{L^2(Q)}^2 + \norma{h^0}_V^2 \big).}
\end{align*}
Therefore, using parabolic regularity and the fact that
\begin{equation}
z(t) =  \frac 1 \a \int_0^t e^{-\beta (t -s) /\alpha}y(s) \pico{\ds} , \quad t\in [0,T],\non
\end{equation}
 we can deduce that 
\begin{align}
	\norma{\pier{z}}_{\H2 H \cap \W{1,\infty} V \cap \H1 W} 
	\leq \pier{\big( \norma{h}_{L^2(Q)}^2 + \norma{h^0}_V^2 \big)}.
\label{pier25}
\end{align}}%

\pier{A combination of the estimates \eqref{pier19}--\eqref{pier25} concludes the proof, since the continuous embedding of $\UU \subset L^2(Q) \times V$, namely, 
\begin{align*}
	\norma{h}_{\L2 H} + \norma{h^0}_V \leq C \norma{\bh}_\UU \quad \hbox{for every } 
	\bh=(h, h^0) \in \UU,	
\end{align*}
ensures} that \eqref{fre:final} is fulfilled.
\end{proof}

\begin{remark}\label{REM:FRE}
{\rm Let us point out that the \Frechet\ differentiability of $\S$ at the fixed control pair $\opt$ is defined from \pier{$\UU_R$} into $\X$
and not from \pier{an open bounded subset of $L^2 (Q) \times \pier V $}, as it may \juerg{appear} (incorrectly) from the estimates above. 
The reason is that for controls $(\bu,\bv)$ just in $L^2(Q)\times \pier V $ 
we \juerg{cannot guarantee} the existence of a strong solution (cf. Theorem~\ref{THM:EX:STRONG}). 
Nevertheless, the above estimates show that, due to the \pier{density of the embedding} of $\UU$ in $L^2(Q) \times \pier{V}$, the \Frechet\ derivative $D\S \opt \in {\cal L}(\UU, \X)$ can be continuously extended to a linear and
\pier{continuous operator from $L^2(Q) \times \pier{V}$ into $\X$.} In particular, denoting that extension with the same symbol $D\S \opt$, the identity $D\S \opt(\bh) = (\xi,\et)$ continues to hold also for \pier{$\bh = (h,h^0) \in L^2(Q) \times \pier{V}$}.}
\end{remark}

\juerg{It is now} a standard matter to derive the first-order \pier{optimality conditions} for \CP\ by combining \eqref{formal:foc}, Theorem \ref{THM:CON:FRE}, and the chain rule. 

\begin{theorem}\label{THM:OPT:FIRST}
Suppose that \ref{ass:1:constants}--\ref{ass:3:pi}, \ref{ass:strong:gamma:1}--\ref{ass:4:pidiff}, \pico{\ref{ass:reg:non:control}--\ref{ass:control:Uad}} are satisfied. 
\pier{Moreover, let the initial data $\ph_0$ and $w_0$ satisfy \eqref{ass:exweak:initialdata}, \eqref{ass:exstrong:initialdata}, \eqref{strong:sep:initialdata}, and let \an{$(\bu, \bv) $ be an optimal control} with $(\bvp,\bw)=\S(\bu, \bv)$.} Then the optimal pair $\opt$ necessarily fulfills the variational inequality 
\begin{align}\label{cn:first}
		\notag
		&  {k_1} \intQ (\bvp - \vp_Q) \xi
		+  {k_2}  \iO (\bvp(T) - \vp_\Omega) \xi(T)
		+  {k_3} \intQ (\bw - w_Q) \et
		+  {k_4}  \iO (\bw(T) - w_\Omega) \et(T)
			\\ & \quad 	\non 	
		+  {k_5} \intQ (\dt \bw - w_Q')  \dt \et
		+  {k_6} \iO (\dt \bw(T) - w'_\Omega) \dt \et(T)
		+  {\nu_1} \intQ \bu(u- \bu) 
		\\ & \quad	
		+  {\nu_2} \iO \Big( {\bv} (v_0 - \bv) + \nabla {\bv} \cdot \nabla (v_0 - \bv) \Big)
		 \geq 0 \quad \forall (u,v_0) \in \Uad,
\end{align}
where $(\xi,\et)$ denotes the unique solution of the linearized system \Lin\ associated \pier{with} the choice $\bh =( u-\bu,v_0 - \bv)$.
\end{theorem}

\pier{We now want to rewrite the optimality conditions in terms of the solution to the adjoint 
problem, in order to simplify the above variational inequality. The backward-in-time system characterizing the adjoint problem} is given, in a strong form, by
\begin{alignat}{2}
	\non
	& - \dt p - \pi(\bvp)\, \dt q 
	- \Delta p  + \gamma'(\bvp)\, p + 
	\an{\tfrac 2 \thc } \pi'(\bvp)\, p
	- \an{\tfrac 1{\thc^2}} \dt\bw \, \pi'(\bvp)\, p &&
 	\\ & \qquad \label{adsys1}
 	= k_1 (\bvp - \vp_Q) \quad&&\mbox{in }\,Q,\\
	& \non
	-\dt q - \a \Delta q \pcol{{}+{}}  \b \Delta (1 \circledast q) - \an{\tfrac 1{\thc^2}} \pi(\bvp)\, p &&
	\\ & \qquad \label{adsys2}
	= k_3 ( 1 \circledast (\bw - w_Q) )
	+ k_5 (\dt \bw - w_Q')
	+ k_4 (\bw(T) - w_\Omega)
	  \quad&&\mbox{in }\,Q,\\
\label{adsys3}
& \dn p= \dn q = 0\quad&&\mbox{on }\,\Sigma,\\[2mm]
\non
& p(T)= {k_2}(\bvp(T) - \vp_\Omega) - {k_6}\pi(\bvp(T)) (\dt \bw(T) - w'_\Omega) , &&
	\\ \label{adsys4}
 & \quad  
 q(T)={k_6} (\dt \bw(T) - w'_\Omega) \quad &&\mbox{in }\,\Omega,
\end{alignat}
\Accorpa\Adj {adsys1} {adsys4}
where \pier{the product $\circledast $ is defined in \eqref{intr2}}.  For convenience, let us \pier{denote by} $f_q$ the source term in \eqref{adsys2}, that is,
\begin{align*}
	f_q:= k_3 ( 1 \circledast (\bw - w_Q) )
	+ k_5 (\dt \bw - w_Q')
	+ k_4 (\bw(T) - w_\Omega) 
\end{align*}
\an{\pico{and notice that the last part $k_4 (\bw(T) - w_\Omega)$} is constant in time.
Moreover, due to \ref{ass:control:target} and to the fact that $\bw$ is a strong solution in the sense of Theorem~\ref{THM:EX:STRONG}, \pier{$f_q$ satisfies}
\begin{align}\label{est:fp}
	\pier{\norma{f_q}_{\L2 H } \leq C (\norma{\bw}_{{\H2 H \cap \W{1,\infty} V \cap \H1 {\Hx2}}}
	 +1 ) \leq C,}
\end{align}
where the above constant certainly depends on $T$.}

The above system reveals why we did also include the possibly redundant objective terms associated to $k_3$ and $k_4$ in \eqref{cost}.
Indeed, the way they appear in the adjoint system above is nonstandard.
\pcol{Another remark concerns the fact that only first-order time derivatives appear in 
\eqref{adsys1}--\eqref{adsys2}, while the corresponding state system, as well as the linearized one, contains an equation with a second-order time derivative as well. However, note that if \eqref{adsys2} is interpreted as an equation in the \juerg{time-integrated} variable $1 \circledast q$, then it turns out that $ - \dt q = \partial_{tt} (1 \circledast q)$, and 
the system~\eqref{adsys1}--\eqref{adsys4} looks more natural}. 

The well-posedness result, as well as the notion of solution to the above system, is specified in the following theorem.

\begin{theorem}\label{THM:ADJ}
Assume that \ref{ass:1:constants}--\ref{ass:3:pi}, \ref{ass:strong:gamma:1}--\ref{ass:4:pidiff}, \pico{\ref{ass:reg:non:control}--\ref{ass:control:target}} hold true. 
\pier{Let the initial data $\ph_0$ and $w_0$ satisfy \eqref{ass:exweak:initialdata}, \eqref{ass:exstrong:initialdata}, \eqref{strong:sep:initialdata}, and 
let $\opt \in \Uad$} be an optimal control for \CP\ with the associated state $(\bvp,\bw)= \S\opt$.  Then the adjoint system \Adj\ admits a unique weak solution $(p,q) $ \pier{with
\begin{align}
	\label{reg:adj:p}
	p & \in \H1 \Vp \cap \L\infty H \cap \L2 V,
	\\
	\label{reg:adj:q}
	q & \in \H1 H \cap \L\infty V \cap \L2 W,
\end{align}
\juerg{that satisfies} the variational equalities 
\begin{align}
	& - \< \dt p, v>_V
	- \iO \pi(\bvp)\, \dt q \, v
	+ \iO \nabla p \cdot \nabla v  
	+ \intQ \gamma'(\bvp) \,p \, v
		\non \\ & \quad 
	+ \frac 2 \thc \intQ \pi'(\bvp)\,p \, v
	- \frac 1{\thc^2} \iO \dt\bw \,\pi'(\bvp)\, p \, v
	= \iO k_1 (\bvp - \vp_Q)\, v,
	\label{pier31}\\
	& - \iO \dt q \, v
	+ \a \iO \nabla  q \cdot \nabla v 
	\pcol{{}-{}} \b  \iO \nabla (1 \circledast q ) \cdot \nabla v
	- \frac 1{\thc^2}  \iO \pi(\bvp) p \, v
= \iO f_q \zeta, \label{pier32}
	\end{align}
for every $v\in V$, almost everywhere in $(0,T)$, and 
the final conditions}
\begin{align}
	p(T) & = {k_2}(\bvp(T) - \vp_\Omega) - {k_6}\pi(\bvp(T)) (\dt \bw(T) - w'_\Omega)  \hspace{-1.5cm}&& \aeO, \label{pier33}
		\\ 
 q(T) & ={k_6} (\dt \bw(T) - w'_\Omega)\label{pier34}
   &&\aeO.
\end{align}
\end{theorem}

\begin{proof}[Proof of Theorem \ref{THM:ADJ}]
We again proceed formally by pointing out the estimates that will imply the existence of a solution. 
These computations can however easily be reproduced in a rigorous framework.
Moreover, before moving on, let us set $Q_t^T := \Omega \times (t,T)$.

\noindent
{\bf First estimate:}
We \pier{take $v=p$ in \eqref{pier31}, $v = - \thc^2 \dt q$ in \eqref{adsys2},}
add the resulting equalities and \pier{note that two terms cancel out. Then,  
\juerg{integration} over $(t,T)$ and by parts yields}
\begin{align}
	& \non
	 \frac 12\IO2 p
	+ \Qtt |\nabla p|^2
	+ \Qtt \gamma'(\bvp) |p|^2
	+ \thc^2 \Qtt |\dt q|^2	
	+\frac {\a \thc^2}{2} \IO2 {\nabla q}
	\\ & \quad \non
		= 
	\frac 12\IOT2 p
	+\frac {\a \thc^2}{2} \IOT2 {\nabla q}
	\an{{}+k_1 \Qtt(\bvp - \vp_Q) p
	- \frac 2 {\thc} \Qtt\pi'(\bvp) \, p^2}
	\\ & \qquad 
	\an{+ \frac 1 {\thc^2} \Qtt\dt \bw \,\pi'(\bvp) \, p^2
	+ \pier{\b \thc^2 \Qtt \nabla (1 \circledast q ) \cdot \nabla (\dt q)}
	- \thc^2 \Qtt f_q \, \dt q.}
	\label{adj:est}
\end{align}
\pier{Notice that the third term on the \lhs\ is nonnegative due to 
the monotonicity of $\gamma$. 
As for the \pico{sixth} term on the \rhs, \pico{we note that $(1 \circledast q )(T)=0 $
in $\Omega$, thus the Young and \Holder\ inequalities allow us to deduce} that
\begin{align*}
		&  \pier{{} \b \thc^2 \Qtt \nabla (1 \circledast q ) \cdot \nabla (\dt q)}
		\\ & \quad 
		= 
		   - \b \thc^2 \iO  \nabla \pier{(1 \circledast q ) (t)}  \cdot \nabla q(t)
		\pier{{}+{}} \b \thc^2 \Qtt |\nabla q|^2
\\ &\quad{}	\leq 	
\frac {\a \thc^2}{4} \IO2 {\nabla q}
			+ C \Qtt  |\nabla q|^2 .	
\end{align*}
\juerg{Concerning} the \pico{third and last terms} on the \rhs, we recall that $(\bvp, \bw)$ 
\juerg{satisfies} \eqref{separation}--\eqref{reg:strong} and that \ref{ass:control:target} and \eqref{est:fp} hold as well. Hence, \juerg{it follows} from Young's inequality that}
\begin{align*}
	 k_1 \Qtt(\bvp - \vp_Q) p
	- \thc^2 \Qtt f_q \, \dt q
	\leq 
	\frac {\thc^2}2 \Qtt |\dt q|^2
	+ C \Qtt \big(|p|^2 +1 \big).
\end{align*}
\pier{Still on the \rhs, the first terms involving the terminal 
conditions are bounded by a constant due to 
\eqref{adsys4} and \ref{ass:control:target}, 
while for the \an{remaining} terms} we owe to the fact that \juerg{ $\bvp, \dt \bw \in L^\infty(Q)$} 
(cf. Theorem~\ref{THM:EX:STRONG}). \pier{Hence, with the help of 
\ref{ass:reg:non:control}, we have} that
\begin{align*}
	- \frac 2 {\thc} \Qtt\pi'(\bvp) \, p^2
	+ \frac 1 {\thc^2} \Qtt\dt \bw \,\pi'(\bvp) \, p^2
	\leq C \Qtt |p|^2.
\end{align*}
Upon collecting the above computations, we \pier{can apply the Gronwall lemma and infer} that 
\begin{align*}
	\norma{p}_{\L\infty H \cap \L2 V}
	+\norma{q}_{\H1 H \cap \L\infty V}
	\leq C.
\end{align*}

\noindent
{\bf Second estimate:}
Next, we proceed with comparison in equation \eqref{adsys2} to deduce that
\begin{align*}
	\pier{\big\|\Delta \big(\a q + \b (1 \circledast q) \big)\big\|_{\L2 H }} \leq C.
\end{align*}
Then, setting $g = \a q + \b (1\circledast q)$, the elliptic regularity theory entails that $\norma{g}_{\L2 W} \leq C .$ Hence, solving the equation $\a q + \b (1\circledast q)=g$ with respect to $1\circledast q$ (which is equal to $0$ at the time $T$), 
we eventually obtain~that
\begin{align*}
	\norma{1 \circledast q }_{\L2 W} + \norma{q }_{\L2 W} \leq C. 
	\end{align*}
\noindent
{\bf Third estimate:}
\juerg{Finally}, we take an arbitrary test function $v \in \L2 V$ \pier{in \eqref{pier31} and compare the terms.} Using the above estimates, it is then a standard matter \juerg{to realize} that
\begin{align*}
	\norma{\dt p}_{\L2 \Vp} \leq C.
\end{align*}

This concludes the proof. In fact, let us recall that the above estimates also imply the uniqueness of the weak solution, as the system \eqref{pier31}--\eqref{pier34} is linear.
\end{proof}

By combining Theorem \ref{THM:OPT:FIRST} with Theorem \ref{THM:ADJ}, we can obtain a more effective version of the variational inequality \eqref{cn:first}.
\begin{theorem}\label{THM:OPT:FINAL}
Suppose that \ref{ass:1:constants}--\ref{ass:3:pi}, \ref{ass:strong:gamma:1}--\ref{ass:4:pidiff}, and \ref{ass:reg:non:control}--\ref{ass:control:Uad} are satisfied. 
\pier{Moreover, assume that the initial data $\ph_0$ and $w_0$ satisfy \eqref{ass:exweak:initialdata}, \eqref{ass:exstrong:initialdata}, \eqref{strong:sep:initialdata}, and 
let $\opt \in \Uad$ be an optimal control for \CP\ with  associated state $(\bvp,\bw)= \S\opt$. Finally, let} $(p,q)$ be the unique solution to the adjoint system \Adj\ as given 
by Theorem~\ref{THM:ADJ}.
Then the optimal pair $\opt$ necessarily verifies 
\begin{align}\label{cn:final}
	 \non &\intQ ( q + {\nu_1} \bu )(u- \bu) 
		+   \iO (q(0) + {\nu_2}{\bv} )(v_0 - \bv) 
		\\ & \quad + {\nu_2} \iO \nabla {\bv} \cdot \nabla (v_0 - \bv)
		 \geq 0 \quad \forall (u,v_0) \in \Uad.
\end{align}
\end{theorem}

\begin{remark}
{\rm Let us point out that the regularity in \eqref{reg:adj:q} entails that $q \in \C0 H$,
 so that $q(0)$ makes sense in $\Lx2.$}
\end{remark}

\begin{proof}[Proof of Theorem \ref{THM:OPT:FINAL}]
Starting from Theorem \ref{THM:OPT:FIRST} \pier{and comparing \eqref{cn:first} with \eqref{cn:final}, we realize that, in order to prove Theorem~\ref{THM:OPT:FINAL}, it  suffices to check that}
\begin{align}\label{claim}
	\non  \intQ qh \pier{{}+{}} \iO q(0) h^0 
	 &{}\pier{{}\geq{}}
	{k_1} \intQ (\bvp - \vp_Q) \xi
		+  {k_2}  \iO (\bvp(T) - \vp_\Omega) \xi(T)
		\\ &\quad{}		\non
		+  {k_3} \intQ ( \bw - w_Q)   \et
		+  {k_4} \iO ( \bw(T) - w_\Omega)  \et(T)
		\\ & \quad {}		
		+  {k_5} \intQ (\dt \bw - w_Q')  \dt \et
		+  {k_6} \iO (\dt \bw(T) - w'_\Omega) \dt \et(T),
\end{align}
with $(\xi,\eta)$ \pico{denoting} the unique solution to \Lin\ associated with the increment $(h,h^0) = (u - \bu,v_0 -\bv)$.
\juerg{To this end}, we test \eqref{linsys1} by $p$, \eqref{linsys2} by $q$, \juerg{and}
integrate over time and by parts to infer that
\pier{%
\begin{align*}
	0 &= \intQ  \Big[ \dt\xi-\Delta\xi + \gamma'(\bvp)\xi + \tfrac2 {\thc} \pi'(\bvp)\xi 
	- \tfrac1{\thc^2}{\dt \et}\,\pi(\bvp)
	- \tfrac1{\thc^2}{\bdtw}\,\pi'(\bvp)\xi\Big] p
	\\ & \qquad  +\intQ  \Big[ \ddt \et -\alpha \Delta (\dt \et)-\beta \Delta \et + \pi' (\bvp)\xi \dt \bvp+ \pi (\bvp) \dt
	 \xi\Big] q
	- \intQ hq
	\\ &  = 
	{}- \int_0^T \!\< \dt p, \xi>_V \juerg{{\rm d}t} + \intQ\! \nabla p \cdot \nabla \xi 
	\\ & \qquad 
	+{}
	\intQ\!  \Big[\gamma'(\bvp) p + \tfrac2 \thc \pi'(\bvp) p
	- \tfrac1{\thc^2} \dt\bw \, \pi'(\bvp)\, p 
	- \dt q \, \pi(\bvp)
 	  \Big] \xi 
	\\ & \qquad 
	+ \intQ\Big[ -\dt q \,  \dt \et + \a \nabla q \cdot \nabla ( \dt \et ) \pcol{{}+{}} \b \nabla (1 \circledast q) \cdot \nabla ( \dt \et )- \tfrac1{\thc^2} \pi(\bvp)  p \, \dt\et \Big]
	\\ &  \qquad 
	+ \iO \Big[ p(T)\xi(T) + \dt \et (T) q(T) + \pi(\bvp(T))\xi(T) q(T)\Big]
	\\ &  \qquad 
	- \intQ qh - \iO q(0) h^0.
\end{align*}
By using \pier{\eqref{pier31}--\eqref{pier34}}, we simplify the above identity, obtaining that}
\begin{align*}
	0&=  k_1\intQ   (\bvp - \vp_Q) \xi
		+ {k_2}\iO (\bvp(T) - \vp_\Omega) \xi(T) 
	\\ & \quad 
	+ \intQ \Big(k_3 ( 1 \circledast (\bw - w_Q) )
	+ k_4 (\bw(T) - w_\Omega)\Big)\dt \et
	\\ &  \quad 
	+ 	k_5 \intQ (\dt \bw - w_Q') \dt \et
	+ 	k_6 \iO (\dt \bw(T) - w_\Omega') \dt \et(T)
	\\ &  \quad 
	- \intQ qh - \iO q(0) h^0. 
\end{align*}
\pier{Now, we integrate by parts the second line, using the initial condition $\et(0)=0$ and the fact that $1 \circledast (\bw - w_Q) (T) =0$.  Then, it is shown that \eqref{claim} holds, and the \pier{proof} is concluded.}
\end{proof}

Finally, let us notice that from \juerg{\eqref{cn:final} we obtain  the standard characterization
 for the minimizers $\bu$ and $\bv$ if $\nu_1$ and $\nu_2$ are positive.} \pier{Prior to the statement, we recall the definition \eqref{Uad} of $\Uad$.}
\begin{corollary}
\label{COR:OP:POINTWISE}
Suppose that the assumptions of Theorem~\ref{THM:OPT:FINAL} hold, and let $\nu_1 > 0 $.
Then, $\bu$ is the $\L2 H$-orthogonal projection of $-{\nu}^{-1}_1 q$ onto the closed and convex subspace $\{ u \in L^\infty(Q) : u_* \leq u \leq u^* \,\,\aeQ\}$, and 
\begin{align}
	\non
	\bu(x,t)=\max \big\{ 
	u_*(x,t), \min\{u^*(x,t),-{\nu}_1^{-1} q(x,t)\} 
	\big\} \quad \hbox{for $a.a. \,(x,t) \in Q.$}
\end{align}
Likewise, if $\nu_2 >0 $, \juerg{then we} infer from Stampacchia's theorem (see, e.g., \cite[Thm.~5.6, p.~138]{Bre}) that $\bv$ is characterized \juerg{by 
}
\begin{align*}
	\frac{\nu_2}2 \norma{\bv}^2_V + \iO q(0)\bv  = \min_{v_0 \in {\mathcal C}} \Big\{ \frac{\nu_2}2 \norma{v_0}^2_V   + \iO q(0)v_0   \Big\},
\end{align*}
where ${\mathcal C} $ denotes the nonempty, closed and convex subset 
$$\{ v_0 \in V : \ \, v_* \leq u \leq v^* \,\,\aeO, \ \, \norma{v_0}_V \leq M\}.$$
\end{corollary}

\section*{Acknowledgments}
This research was supported by the Italian Ministry of Education, 
University and Research~(MIUR): Dipartimenti di Eccellenza Program (2018--2022) 
-- Dept.~of Mathematics ``F.~Casorati'', University of Pavia. 
In addition, \pcol{PC and AS gratefully mention their affiliation}
to the GNAMPA (Gruppo Nazionale per l'Analisi Matematica, 
la Probabilit\`a e le loro Applicazioni) of INdAM (Isti\-tuto 
Nazionale di Alta Matematica). \pcol{Moreover, PC aims to point out} his collaboration,
as Research Associate, to the IMATI -- C.N.R. Pavia, Italy.

\End{document}
